\documentstyle[12pt,epsfig]{article}
\newcommand{\bt}{\begin{Theorem}}
\newcommand{\et}{\end{Theorem}}

\newcommand{\ei}{\end{itemize}}
\newcommand{\bea}{\begin{eqnarray}}
\newcommand{\eea}{\end{eqnarray}}
\newtheorem{Theorem}{\sc Theorem}
\newtheorem{Lemma}[Theorem]{\sc Lemma}
\newtheorem{Proposition}[Theorem]{\sc Proposition}
\newtheorem{Corollary}[Theorem]{\sc Corollary}
\newtheorem{Definition}[Theorem]{\sc Definition}
\newtheorem{Example}[Theorem]{\sc Example}
\newtheorem{Remark}[Theorem]{\sc Remark}

\newcommand{\be}{\begin{equation}}
\newcommand{\ee}{\end{equation}}

\def\qed{\hfill$\Box$}
\def\circleds{{\bigcirc \!\!\!\!s}}

                     \def\R{{I\!\!R}}

                     \def\CZ{{\rm\ke}rn.26em
                     \newcommand\la{{\langle}}
                     \newcommand\ra{{\rangle}}
                     \newcommand\lar{\leftarrow}
                     \newcommand\Lar{\Leftarrow}
                     \newcommand\rar{\rightarrow}
                     \newcommand\Rar{\Rightarrow}

                     \vrule width.02em
                     height.5ex depth0ex
                     \kern.04em
                     \vrule width .02em
                     height1.47ex depth-1ex
                     \kern-.34em Z}
                     
                     \def\C{{\rm \kern.24em
                     \vrule width.02em
                        height1.4ex depth-.05ex
                     \kern-.26em C}}

                     \def\ra{{\rightarrow}}

                     \def\ei{{\bf e_i}}

                     \def\N{{\rm I\kern-.23em N}}
                     \def\B{{\rm I\kern-.25em B}}
                     \def\D{{\rm I\kern-.25em D}}
                     \def\E{{\rm I\kern-.25em E}}
                     \def\F{{\rm I\kern-.25em F}}

                     \def\I{{\rm I\kern-.25em I}}
                     
                     \def\M{{\rm I\kern-.23em M}}
                     \def\P{{\rm I\kern-.25em P}}
                     \def\A{{\rm \kern.22em
                     \vrule width.02em
                        height0.5ex depth 0ex
                     \kern-.24em A}}
                     \def\G{{\rm \kern.24em
                     \vrule width.02em
                        height1.4ex depth-.05ex
                     \kern-.26em G}}
                     \def\J{{\rm \kern.19em
                     \vrule width.02em
                        height1.47ex depth 0ex
                     \kern-.21em J}}
                     \def\O{{\rm \kern.24em
                     \vrule width.02em
                        height1.4ex depth-0.5ex
                     \kern-.26em O}}
                     \def\Q{{\rm \kern.24em
                     \vrule width.02em
                        height1.4ex depth-.05ex
                     \kern-.26em Q}}
                     \def\S{{\rm \kern.18em
                     \vrule width.02em
                        height1.4ex depth-.9ex
                     \kern.12em
                     \vrule width.02em
                     height0.7ex depth 0ex
                      \kern-.34em S}}
                     \def\T{{\rm \kern.45em
                     \vrule width.02em
                        height1.47ex depth 0ex
                     \kern-.47em T}}
                     \def\U{{\rm \kern.30em
                     \vrule width.02em
                        height1.47ex depth-.05ex
                     \kern-.32em U}}
                     \def\V{{\rm \kern.27em
                     \vrule width.02em
                        height1.47ex depth-.8ex
                     \kern-.29em V}}
                     \def\W{{\rm \kern.25em
                     \vrule width.02em
                        height1.47ex depth-0.9ex
                     \kern.34em
                     \vrule width.02em
                        height1.47ex depth-.9ex
                      \kern-.63em W}}
                     \def\X{{\rm \kern.30em
                     \vrule width.02em
                        height1.4ex depth-1ex
                      \kern.12em
                      \vrule width.02em
                         height0.4ex depth 0ex
                      \kern-.46em X}}
                     \def\Y{{\rm \kern.25em
                     \vrule width.02em
                        height1.0ex depth 0ex
                     \kern-.27em Y}}
                     \def\Z{{\rm \kern.26em
                     \vrule width.02em
                        height0.5ex depth 0ex
                      \kern.04em
                      \vrule width.02em
                         height1.47ex depth-1ex
                      \kern-.34em Z}}

\begin{document}

\begin{center}
{\large{\bf On Product Systems arising from  Sum Systems}}\\

\vskip0.5em

{\em B. V. RAJARAMA BHAT},

Indian Statistical Institute,\\
Bangalore-560 059, India.\\
bhat@isibang.ac.in

{\em R. SRINIVASAN }\\
Department of Mathematical Sciences,\\
University of Tokyo, Komaba, Tokyo, 153-8914, Japan.\\
vasanth@ms.u-tokyo.ac.jp

\today
\end{center}

\begin{abstract}
Boris Tsirelson constructed an uncountable family of type $III$ product
systems  of Hilbert spaces 
through the theory of Gausian spaces, measure type spaces and
`slightly
coloured noises', using  techniques from probability theory. Here we
take a purely functional analytic approach and try to have a better
understanding of Tsireleson's construction and his  examples.

We prove an extension of Shale's theorem connecting symplectic
group and Weyl representation. We show that the  `Shale map'
respects compositions (This settles an old conjecture of
K. R. Parthasarathy \cite{KRP}).
 Using this we associate
a product system to a sum system. This construction includes the
exponential product system of Arveson, as a trivial case, and the
type $III$ examples of Tsirelson. 

By associating a von Neumann algebra to
every `elementary set' in $[0,1]$, in a much simpler and direct
way, we arrive at the invariants of the product system introduced by
Tsirelson,
given in terms of
the sum system. Then we introduce a notion
of
divisibility for a sum system, and prove that the examples of
Tsirelson are divisible. It is shown that only type $I$ and
type $III$ product systems arise out of divisible sum systems. 
Finally, we give a sufficient condition for
a divisible sum system to give rise to a unitless (type $III$)
product system.
\end{abstract}

\bigskip \noindent
\section{Introduction:}

R. T. Powers \cite{Po1} initiated a study of $E_0-$semigroups, which
are weakly continuous semigroups of unital 
$*-$endomorphisms of some $\B(H)$,
for a separable Hilbert space $H$. In this context 
Arveson \cite{Arv} introduced the concept of
product system of Hilbert spaces as an invariant for
$E_0$-semigroups. Up to cocycle conjugacy an
$E_0$-semigroup $\{\alpha_t\}$ is determined by the family of Hilbert
spaces $\{H_t\}$,
where $$H_t =\{T \in \B(H):\alpha_t(X) T = TX, ~ \forall X \in
\B(H)\}$$ with inner product $\langle T, S \rangle 1_H= S^*T$ (see
\cite{Arv}). Moreover the family $\{H_t\}$ forms a
product
system of
Hilbert spaces (see Definition \ref{productsystem}). Arveson also
constructed an $E_0-$ semigroup from a given product system, thus
proving
that the product systems forms a complete invariant for the
$E_0-$semigroup (up to cocycle conjugacy).

Arveson classified product systems, according to the
existence of
units (see Definition \ref{unit}), into three broad categories,
such as
type
$I$, $II$, $III$. He also classified completely the type $I$
product systems,
up to
isomorphism. We refer to \cite{Arv3} for general theory of
$E$-semigroups and product systems and \cite{Pr} for some
recent developments.

The theory of product systems
was lacking  enough examples. For quite sometime there were 
essentially only one
example  each for type $II$ and type $III$ product systems (due to
R. T. Powers (see \cite{Po1}-\cite{Po3})). 
Tsirelson produced an uncountable family of both type $II$ and type
$III$
product systems (ref \cite{Ts1}, \cite{Ts}).

Tsirelson uses the theory of random sets arising from a Brownian motion
to get type $II$ product systems and the theory of FHS spaces, Gaussian
spaces,
measure type spaces and what he calls as `slightly coloured noises' to
get the examples of type $III$ product systems. Tsirelson's
construction
of type $III$ product systems is complicated and involves lots of techniques from
probability theory. Also it is not clear as how to
work with the $E_0-$semigroup associated with the product systems, and
there is no information regarding other invariants of the
product system, such as the automorphism group etc.
Our work is inspired by the path
breaking results of Tsirelson (which in turn borrow on
some brilliant ideas of Vershik).

 The basic idea of Tsirelson's construction of
type $III$ product systems is simple. Usual $L^2$ on sub-intervals on real
line is a direct sum system in the sense that $L^2(0, s)\oplus L^2(s, s+t)=
L^2(0, s+t)$ for positive $s, t.$ Such a system on
 `exponentiation' gives the type I or 
the Fock product system. Now if we replace a direct sum system by an
`almost' or `quasi' direct sum system we get more exotic product systems.
First job is to make precise as to what one means by quasi-direct sum and
then one has to find a suitable procedure of exponentian. Tsirelson does
this by his notion of FHS equivalence, identifying the Hilbert spaces in the
sum system with Gaussian type spaces, and then getting 
the product system, as
the
$L^2$-space of the corresponding measure type spaces.
 We retain Tsirelson's notion of sum system
though we don't use the language of probability theory.
The essential difference in our approach is that 
we do the exponentian using the theory of symmetric Fock
spaces and
a generalised version of Shale's theorem.

 We
first prove
in Section \ref{sumprod}, a generalisation of Shale's theorem.
 We also
prove a
functorial property in the Shale's theorem affirmatively settling
a conjecture of K. R. Parthasarathy.
 Using this, after proving some lemmas, we
associate product system with a sum system. We show that this gives the
exponential product system, as a trivial case, and includes the examples
of Tsirelson. We also prove some properties of sum systems, and provide
an operator theoretic proof of some facts in Tsirelson's work.

In Section \ref{non-isomorphic}, given a product system 
we associate a von Neumann algebra to
any elementary set (finite union of intervals) in the interval $[0,1]$.
  We analyse these von Neumann algebras, and by simple
application of double commutant theorem, strong- weak convergences, we
arrive at the invariants
for the
product systems, given in terms of the original sum systems.  This would
prove the examples
of Tsirelson are non-isomorphic to each other. This is infact the
difficult part of Tsirelson's work, and we give here a much more direct
and simple
proof of this fact.

In Section \ref{onunits}, we first define a notion of divisibility for a
sum system, and study some basic properties of a divisible sum system. We
prove that all examples of Tsirelson are divisible. We also show that
only type $I$ and $III$ are possible under the divisibility assumption on
the sum system. Finally, using some of the notions introduced by
Tsirelson, we prove a sufficient condition for the product system
arising from a divisible sum system to be
of type $III$.

Almost after finishing this work we came to know about the new
preprint of
Tsirelson (\cite{Ts3}), where he has simplified many of the proofs in his
earlier two preprints, for producing the uncountable family of type $II$
and
$III$ product systems. We still believe that our method is more direct and
simple, leading to new applications. We plan to consider some research ideas
emerging from this approach in future.

We end this section by recalling some of basic definitions, which are
intially defined by Arveson. For undelying measurability conditions
 we use a slightly modified, but
essentially equivalent, definition, given by Volkmar Liebcher
(\cite{Vol}).

\begin{Definition}\label{productsystem}
A product system of Hilbert spaces is an one parameter family
of separable Hilbert spaces $\{H_t\}_{t \in (0,\infty)}$, together with
unitary
operators $$U_{s,t} : H_s \otimes
H_t ~\mapsto H_{s+t}~ \mbox{for}~ s, t \in (0,\infty),$$
satisfying the following two axioms of associativity and measurability.

\medskip
\noindent(i) (Associativity) For any $s_1, s_2, s_3 \in (0,\infty)$
$$U_{s_1, s_2 + s_3}( 1_{H_{s_1}} \otimes U_{s_2 ,
s_3})=  U_{s_1+ s_2 , s_3}( U_{s_1 ,
s_2} \otimes 1_{H_{s_3}}).$$

\medskip
\noindent (ii) (Measurability) There exists a countable set $H^0$ of
sections $ R\ni t \rightarrow h_t \in H_t$ such that $ t  \mapsto 
\langle h_t,
h_t^\prime\rangle$ is
measurable for any two $h, h^\prime \in H^0$, and the set $\{h_t: h \in
H^0\}$ is a total set in $H_t$, for each $ t \in (0,\infty)$.
Further it is also assumed that the map $(s,t) \mapsto \langle
U_{s,t}(h_s
\otimes h_t), h^\prime_{s+t} \rangle$ is measurable for any two $h ,
h^\prime \in H^0$.
\end{Definition}

\begin{Definition}
Two product systems $(H_t, U_{s,t})$ and $(H^\prime_t, U^\prime_{s,t})$
are
said to be isomorphic if there exists a unitary operator $V_t:H_t \mapsto
H_t^\prime$, for each
$t \in (0,\infty)$, satisfying the following two conditions.

\medskip\noindent
(i) $V_{s+t}U_{s,t}= U_{s,t}^\prime (V_s \otimes V_t).$

\medskip\noindent
(ii) The $t \in (0,\infty) \mapsto
\langle V_t h_t,
h^\prime\rangle$ is measurable for any $ h \in H^0, h^\prime \in
{H^\prime}^0$.
\end{Definition}

\begin{Remark}
Volkmar Liebcher has proved in \cite{Vol} that any two
measurable structures give rise to isomorphic product systems,
and as a consequence we get that two product systems are isomorphic if they are
algebraically isomorphic. That is the condition (ii) in the above definition
can be dropped.
\end{Remark}

\begin{Definition}
For a product system $(H_t, U_{s,t})$, we define the opposite product
system $(H^{op}_t, U^{op}_{s,t})$ by, $$H^{op}_t = H_t, ~~
U_{s,t}^{op}=U_{t,s}\tau_{s,t},$$ where $\tau_{s,t}$ is the flip operator
on $H_s \otimes H_t$, $\tau_{s,t}(x \otimes y) = y \otimes x.$

A product system is said to be symmetric if it is isomorphic to its
opposite product system, (i.e) it is anti-isomorphic to itself.
\end{Definition}

We next define the units, based on whose existence, the product systems
are classified into three broad categories.

\begin{Definition}\label{unit}
A unit is a measurable section $\{ u_t\}_{t \in
(0,\infty)}$, ((i.e) $u_t \in H_t$, and the map $t \mapsto \langle u_t,
h_t\rangle$ is
measurable for any $h \in H^0$), satisfying
$$U_{s,t}(u_s \otimes u_t)= u_{s+t}, ~  \forall s,t \in (0,\infty),~
\mbox{and} ~u_t \neq 0~ \mbox{for some} ~ t \in (0, \infty).$$
\end{Definition}

\bigskip

We denote by $\U$ the set of all units for a product system. We say a
product system is of type $I$, if units exists for the product
system and they generate the product system, (i.e.) for any fixed $t \in
(0,\infty)$, the set $$\{u^1_{t_1}u^2_{t_2}  \cdots u^n_{t_n}:
\sum_{i=1}^n
t_i = t, u^i \in \U\},$$ is a total set in $H_t$, where the product is
defined as the image of
$u^1_{t_1}\otimes u^2_{t_2}  \cdots \otimes u^n_{t_n}$ in
$H_t$, under
the canonical unitary given by the associativity axiom. It is of type
$II$ if
units exists but they don't generate the product system. We say a product
system to be of type $III$ or unitless if there does not exist any unit
for the product system. We are most concerned about this type $III$
product systems in this paper.

\section{The construction}\label{sumprod}
 In this section we construct a product system from a given sum
system (see definition \ref{sumsystem}). We do this by
proving a  generalised version of Shale's theorem. Before
that we fix our notation.

For a real Hilbert space $G$ we denote by $\overline{G}$ the
complexification of $G$. (Throughout this paper we always denote a real
Hilbert space by $G$, and if the Hilbert space is complex
we denote it by $ H$ or $\overline{G}$ or we specify it). 
We define, for a single Hilbert
space
$G$ or
for two
Hilbert spaces
$G_1$ and
$G_2$, $\S(G)$ and $ \S(G_1, G_2)$ in the following way,
$$\S(G)=\{A\in \B(G): A~ \mbox{positive,
invertible and }~ I-A ~\mbox{is Hilbert-Schmidt} \},$$
$$\S(G_1,G_2)=\{A\in \B(G_1, G_2):
A~\mbox{invertible and} ~I-(A^*A)^{\frac{1}{2}} ~\mbox{Hilbert-Schmidt}
\}.$$ In the above definition, and elsewhere in this paper, by
invertible we mean the inverse is also
bounded.
Note that $\S(G, G)$ is different from $\S(G)$.

Clearly
$A\in \S(G_1,G_2)$ if and only
if $A^{-1} \in \S(G_2, G_1)$, and $A\in
\S(G_1,G_2)$ if and only
${(A^*A)}^{\frac{1}{2}} \in \S(G_1)$.

If $A\in \S(G_1, G_2)$ and $B\in \S(G_2, G_3)$, we may conclude from the
relation
$I-A^*B^*BA=I-A^*A+A^*(I-B^*B)A$ that $I- A^*B^*BA$ is a
Hilbert-Schmidt operator. The fact that $I -A^*A$ is a Hilbert-Schmidt
operator is equivalent to saying that
$I-(A^*A)^{\frac{1}{2}}$ is a Hilbert-Schmidt operator, when $A$ is
invertible (see
\cite{FHS}, and \cite{Ts} Proposition 9.9, page 46), and the above
verification now proves  that $BA \in \S(G_1, G_3)$.

Also if $A \in  \S(G_1, G_2)$, then the same fact
implies that $A^*A \in \S(G_1)$. Now, as $A^{-1} \in \S(G_2, G_1)$ and
$A^*A \in \S(G_1)$,
we conclude that $(A^*)^{-1} \in \S(G_1, G_2)$.

Suppose $A \in \S(G_1, G_2)$, and $G_1^\prime \subset G_1$ be any
subspace and $G_2^\prime
= A(G_1^\prime)$, then we want to check whether the
restricton
$A|_{G_1^\prime} \in \S(G_1^\prime, G_2^\prime)$. 
As $I_{G_1^\prime} - (A|_{G_1^\prime})^* A|_{G_1^\prime}$ is
the compression of $I_{G}-A^*A$ to $G_1^\prime $ it is 
a Hilbert-Schmidt
operator.  Then it is clear  that
$A|_{G_1^\prime} \in \S(G_1^\prime, G_2^\prime)$.

\bigskip

We make some definitions and fix some notation.

\begin{Definition}
We say two subspaces $G_1$ and $G_2$, both contained in a real
Hilbert space $G$, are quasi-orthogonal if there exists a map $A
\in \S(G)$ such that $\langle Ax, Ay\rangle=0, ~ \forall~ x \in G_1, y
\in G_2.$
\end{Definition}

We use the notation
$G_1 \uplus G_2 = G$ (respectively
$\uplus _{i=1}^nG_i =G$), if $G_1$ is
quasi-orthogonal
to $G_2$, and $G$ is
generated by
$G_1$ and $G_2$ ( respectively $G_i$'s are mutually 
quasi-orthogonal and $G= \vee_{i=1}^n G_i$). We also denote by
$$\O(\oplus_{i=1}^n G_i, G)=\{A \in \S(\oplus_{i=1}^nG_i, G); A(G_i)=G_i,
~\mbox{for
each}~i=1, 2, \cdots n\}. $$

\begin{Lemma}\label{ass0}
Let $\{G_i\}_{i=1}^n$ be a family of real
Hilbert spaces all contained in a one
real Hilbert space
$G$. Then the set $\O(\oplus_{i=1}^n \G_i, G)$ is not empty
if and only if $\uplus _{i=1}^n G_i=G$.
\end{Lemma}

{\em Proof:}
Suppose there exists $U \in \O(\oplus^n_{i=1}
\G_i, G)$
then define
$A=$ $((U^{-1})^*U^{-1})^{\frac{1}{2}}$. Then
for $x\in G_i, y\in G_j$
and $i\neq j$, $\langle Ax, Ay\rangle = \langle U^{-1}x,
U^{-1}y\rangle=0$. The invertibility of $U$
and the
condition $U(G_i)=G_i$
                     clearly imply that
$G=span\overline{[\cup_{i=1}^n G_i]}$.

Now to prove the otherway, suppose there
exists $A \in \S(G)$
such that  $AG_i \perp
AG_j$ if $i \neq j$, then  as
$G=span\overline{[\cup_{i=1}^n
G_i]}$ we conclude that $\oplus^n_{i=1} AG_i
=G$. Now define $U=A^{-1}(\oplus^n_{i=1}
(A|_{G_i}))$, clearly $U \in \O(\oplus_{i=1}^n \G_i, G).$
\qed

\begin{Remark}\label{indpt} Note that we have also proved that
 $\O(\oplus_{i=1}^n \G_i, G)$ is not empty if and only if the map
$\oplus_{i=1}^nx_i~ \mapsto ~ \sum_{i=1}^nx_i$ is in $\S(\oplus_{i=1}^n
\G_i, G)$.
\end{Remark}

\bigskip

The following lemma is proved in \cite{Ts} using probability theory
(Radon-Nikodym
derivatives). We provide an operator theoretic proof here.

\begin{Lemma}\label{ass}
Let $G_1, G_2, G_3$ be real Hilbert spaces
all contained in a real Hilbert
space $G$. Let $G_{12}$ (resptly. $G_{23}$)
be the Hilbert space generated
by $G_1$ and $G_2$ (resptly. by $G_2$ and
$G_3$). Suppose that $G_1  \uplus G_{23} =G$
and $G_{12} \uplus G_3 = G$,
then it also holds that $\uplus _{i=1}^3G_i = G$.
\end{Lemma}

{\em Proof:} Choose $A_1$, $A_2 \in \S(G)$ such
that
$A_1G_1 \perp A_1G_{23}$ and $A_2 G_{12}\perp A_2 G_3.$
As we also have
$G=span\overline{[G_{12}, G_3]}$, we conclude
that $A_2 G_{12} \oplus
                     A_2G_3 =G$. Now define $$
A_0 ~=~\left(((A_1{A_2}^{-1}|_{A_2
                     G_{12}})^*A_1{A_2}^{-1}|_{A_2
G_{12}})^{\frac{1}{2}}
\oplus
I|_{A_2G_3}\right) A_2,$$ where
$((A_1{A_2}^{-1}|_{A_2
G_{12}})^*A_1{A_2}^{-1}|_{A_2
G_{12}})^{\frac{1}{2}}
\oplus I$
is defined on $A_2 G_{12} \oplus
A_2G_3$. Clearly $A_0 \in \S(G,G)$.

Now for $x \in G_1, y\in G_2$, we have
$$\langle A_0x,
A_0 y\rangle =\langle (A_1{A_2}^{-1}|_{A_2
G_{12}})^*A_1{A_2}^{-1}|_{A_2G_{12}}
(A_2x),A_2y\rangle_{A_2 G_{12}}$$
$$=\langle A_1 A_2^{-1} (A_2 x), A_1
A_2^{-1}(A_2y)\rangle =\langle A_1
x,A_1 y\rangle =0.$$ Also if $x \in G_{12}$ and $y
\in G_3$ $\langle A_0x, A_0y\rangle =\langle z,
A_2y\rangle=0,$ where $z$ is some element in
$A_2 G_{12}$. So $A_0$ satisfies
$\langle A_0x, A_0y\rangle =0$ whenever $x \in G_i, y
\in G_j$,
for $1\leq i,j \leq 3$ and
$1\neq j$. If we define $A=(A_0^*A_0)^{\frac{1}{2}}$, then clearly $A \in
\S(G)$ and it continues to satisfy $\langle Ax, Ay\rangle =0$ whenever $x
\in G_i, y
\in G_j$,
for $1\leq i,j \leq 3$ and
$1\neq j$.
\qed

 Let $G_1, G_2$ be two real Hilbert spaces and let $A
\in
\S(G_1,G_2 )$, then define
$S_A:\overline{G_1} \rightarrow
\overline{G_2}$ by
$S_A(u+iv)=Au+i(A^{-1})^*v$ for $u,v \in
G_1$. Then $S_A$ is a symplectic
isomorphism between $\overline{G_1}$ and
$\overline{G_2}$
(i.e. $S_A$ is a real linear, bounded,  invertible map with a
bounded inverse satisfying $Im(\langle S_Ax,S_Ay\rangle
)=Im\langle x,y\rangle $ for all $x,y
\in \overline{G_1}$, see
\cite{KRP} page 162). Notice that, for a unitary operator $U \in \B(G_1,
G_2)$ (which is clearly in $\S(G_1, G_2)$), $S_U$ is a complex linear,
unitary operator, and $S_U(x +iy) = Ux + iUy$.

We briefly recall the notions of the symmetric Fock space of a Hilbert
space, exponential vectors
and the Weyl operators. For a complex Hilbert space $K$, we know that the
tensor product $\otimes_{i=1}^n K_i$, where $K_i=K$ for all $i=1,2,
\cdots n$, admits an action of the symmetric group $S_n$, given by
$$\sigma(\otimes \xi_i)= \otimes \xi_{\sigma^{-1}(i)}.$$ The symmetric
tensor
product and symmetric Fock space corresponding to $K$ are defined by
$$K^{\circleds^n}= \{\xi \in K: \sigma(\xi)= \xi\}, ~ \Gamma_s(K)=
\oplus_{i=0}^\infty K^{\circleds^n}, $$ where $K^{\circleds^0}$ is
assumed to be $\C$. We
call $1 \in \C \subset \Gamma_s(K)$, as the vacuum vector, and denote it
by
$\Phi$. For any $x \in K$, we define, $$e(x)= \oplus_{i=0}^\infty
\frac{x^{\otimes^n}}{\sqrt{n!}}.$$ It is a fact that the set $\{e(x): x
\in K\}$ is a
linearly independent and  total set in $\Gamma_s(K)$. The Weyl operator,
corresponding to an element $x \in K$ is defined by, $$W(x)(e(y))=
e^{-\frac{1}{2}\|x\|^2- \langle y, x
\rangle} e(y+x),$$ and $W(x)$ is extends to an unitary operator on
$\Gamma_s(K)$. Also, for a unitary operator $U$, between two Hilbert
spaces
$K_1$ and $K_2$, $U \in \B(K_1, K_2)$, we define another operator
$Exp(U)$ between the corresponding symmetric Fock spaces, $Exp(U) \in
\B(\Gamma_s(K_1), \Gamma_s(K_2))$, by, $$ Exp(U)(e(x))= e(Ux).$$ Again,
$Exp(U)$ extends to an unitary operator.

As $W(x)W(y) = e^{-Im \langle x, y\rangle }W_{x+y}$, the 
correspondence $x \mapsto W(x)$ provides a projective representation
for the abelian group $K$. Notice that
when $K= \overline {G}$, as $Im(\langle S_Ax, S_Ay\rangle )= Im \langle x, y\rangle ,$
the
correspondence $x \mapsto W(S_Ax)$ also provides a projective
representation. Shale's theorem answers the question as to when these two
projective representations are equivalent. 
The following theorem is a
generalisation  of
Shale's Theorem (see \cite{KRP} page 169, Theorem
22.11), where now instead of maps from a real Hilbert space to
itself we have maps from one real Hilbert space to another. 
More importantly we prove that the  `Shale map'
 $\Gamma(\cdot)$, of Shale's
theorem  respects composition - see $(ii)$ of
Theorem \ref{shales}. This was left as an open problem in 
\cite{KRP} page
170). 

\bigskip

\begin{Theorem}\label{shales} (i) Let $G_1,
G_2$ be real Hilbert spaces and $A \in \S(G_1,G_2)$ ,
then there exists a unique unitary operator
                     $\Gamma(A):\Gamma_s(\overline{G_1})
\rightarrow
\Gamma_s(\overline{G_2})$ such that

\begin{eqnarray}\label{SHL}
\Gamma(A)W(u){\Gamma(A)}^* & = & W(S_Au)\label{SHL1}\\
\langle \Gamma(A) \Phi_1, \Phi_2\rangle & \in & \R^{+}\label{SHL2}
\end{eqnarray}
where $\Phi_1$ and $\Phi_2$ are the vacuum
vectors in $\Gamma_s(\overline{G_1})$ and
$\Gamma_s(\overline{G_2})$
respectively.

\medskip \noindent
(ii) Suppose $G_1, G_2, G_3$ be three real Hilbert
spaces, and
$A \in \S(G_1, G_2),~B \in \S(G_2, G_3)$,  then
\begin{eqnarray}\label{shlrep}
\Gamma(A^{-1}) & = & {\Gamma(A)}^*\\\label{shlrep1}
\Gamma(BA) & = & \Gamma(B)\Gamma(A)
\end{eqnarray}

\medskip \noindent
(iii) If $\{T_n\}\subset \S(G, G)$, be any sequence of operators
such that $T_n$
converges strongly to $T \in \S(G, G)$ and $(T_n^*)^{-1}$ converges
strongly to $(T^*)^{-1}$, then $\Gamma(T_n)$ converges
weakly to $\Gamma(T)$.
\end{Theorem}

{\em Proof:} (i) Let $A_0=(A^*A)^{\frac{1}{2}}$. As
$I-A_0$ is a Hilbert-Schmidt operator on
$G_1$, there exists an orthonormal basis $\{e_i\}
\subset G_1$ such that $A_0e_i =\lambda_i e_i$, with $\lambda_i >0$ for
each
$i$ and $\sum_i (\lambda_i -1)^2 < \infty$.

                     Let $f_i= {\lambda_i}^{-1} A e_i$, then as
$$\langle Ae_i, Ae_j\rangle =\langle A^*A
e_i, e_j\rangle =\langle A_0e_i,A_0e_j\rangle = \langle \lambda_i
e_i,\lambda_j e_j\rangle,$$
we conclude that $\{f_i\}$ is an
orthonormal basis for $G_2$.
Also note that $Ae_i=\lambda_i f_i$ and ${A^*}^{-1}e_i ={\lambda_i}^{-1}
f_i$. Now identify
$\overline{G_1}$ with $l^2(\{e_i\})$
((i.e.) with $\oplus_{i=1}^{\infty}\C e_i$)
and
$\overline{G_2}$ with
$l^2(\{f_i\})$ ((i.e) with
$\oplus_{i=1}^{\infty} \C
f_i$). Also
identify $\Gamma_s(\C)$ with $L^2(\R)$ by
$$e(z) \rightarrow
(2\pi)^{-\frac{1}{4}}exp(-\frac{1}{4}t^2+z
t-\frac{1}{2}z^2)$$ for $z \in \C$ (see \cite{KRP}
page 142, Proposition 20.9).
Let $U_i$ (resptly. $V_i$) be the unitary
operator between $\Gamma_s(\C
e_i)$ (resptly. $\Gamma_s(\C f_i)$) and
$L^2(\R)$. Then the following
relations hold (and also with $U_i$
replaced by
$V_i$).
\begin{eqnarray}\label{uni}
(U_i e(z e_i))(t) & =&
(2\pi)^{-\frac{1}{4}}exp(-\frac{1}{4}t^2
+ zt-\frac{1}{2}z^2),\\ \nonumber
(U_i W(x e_i)(U_i)^{-1} f)(t) & =& f(t-2x),\\\nonumber
(U_i W(iy e_i)(U_i)^{-1} f)(t) & =& e^{ity} f(t),
\end{eqnarray}
where $z\in \C$, $f \in L^2(\R)$ and $z=x+iy$
(again see \cite{KRP} page 142,
Proposition 20.9).

                     For $\lambda > 0$ define $L_{\lambda}$ on
$L^2(\R)$ by $L_{\lambda}(f)(x)
= {\lambda}^{- \frac{1}{2}}f(\frac{x}{\lambda})$. $L_{\lambda}$ is a
unitary operator on $L^2(\R)$. Also if we
define $V_{\lambda_i}=
{V_i}^{-1} L_{\lambda_i} U_i$, then clearly
$V_{\lambda_i}$ is a
unitary operator between $\Gamma_s(\C e_i)$
and $\Gamma_s(\C
f_i)$.  Moreover a simple calculation, using
the equations \ref{uni},
shows
that, for any $z=x+iy \in \C$,
$V_{\lambda_i}$
satisfies the following equations.

\begin{eqnarray}\label{weyl}
V_{\lambda_i}W(ze_i){V_{\lambda_i}}^{-1} & =&
W((\lambda_i x
+{\lambda_i}^{-1}y)f_i)\\\label{weyl1}
\langle V_{\lambda_i} \Phi_1,\Phi_2\rangle & = &
({\frac{\lambda_i+{\lambda_i}^{-1}}{2}})^{-\frac{1}{2}},
\end{eqnarray}
where $\Phi_1$ and $\Phi_2$ are the vacuum
vector in $\Gamma_s(\C e_i)$
and $\Gamma_s(\C f_i)$ respectively.

Identify $\Gamma_s(\overline{G_1})$
(resptly.
$\Gamma_s(\overline{G_2})$) with
$\otimes_{i=1}^{\infty}\Gamma_s(\C e_i)$
(resptly. with
$\otimes_{i=1}^{\infty}\Gamma_s(\C f_i)$), where  the
countable tensor product is with respect to
the stabilising sequence of
vacuum vectors. Define $$\Gamma_n= V_{\lambda_1}
\otimes V_{\lambda_2} \otimes \cdots
\V_{\lambda_n} \otimes I_{[n+1},$$
where $I_{[n+1}$ is like an identity operator
between $\otimes_{i=n+1}^{\infty} \Gamma_s(\C e_i)$
and
$\otimes_{i=n+1}^{\infty} \Gamma_s(\C f_i)$,
$$(i.e.)~~ I_{[n+1}
                     (\otimes_{i=n+1}^l z_i e_i \otimes \Phi_1
\otimes \Phi_1 \cdots)~=~
                     \otimes_{i=n+1}^l z_i f_i \otimes \Phi_2
\otimes \Phi_2 \cdots$$
For any $n > m> k$, we have
$$\|(V_{\lambda_{k+1}} \phi_1 \otimes \cdots
\otimes
\V_{\lambda_{n}} \phi_1) \otimes \Phi_2
\otimes \Phi_2 \otimes \cdots
~-~
(V_{\lambda_{k+1}} \phi_1 \otimes \cdots
\otimes
                     \V_{\lambda_{m}} \phi_1) \otimes \Phi_2
\otimes \Phi_2 \otimes
                     \cdots\|^2$$
                     $$~=~2\left(1-\prod_{i=m+1}^n\left(\frac{\lambda_i
+{\lambda_i}^{-1}}{2}\right)^{-\frac{1}{2}}\right),$$
                     which converges to $0$ as $n, m \rightarrow
0$.

For $u \in \overline{G_1}$ define $$\psi(u)=
e^{-\frac{\|u\|^2}{2}}e(u).$$ Clearly $\|\psi(u)\|=1.$ Now we conclude,
 for any $u \in
\oplus_{i=1}^k\Gamma_s(\C e_i)$, that
$$\lim_{n \rightarrow
                     \infty}\Gamma_n(\psi(u))~=~ (V_{\lambda_1}
\otimes V_{\lambda_2}\otimes
\cdots \otimes V_{\lambda_k}\psi(u)) \otimes
\bigotimes_{j=k+1}^{\infty}
                     V_{\lambda_j} \Phi_1$$ exists. Define
$$\Gamma(A)(\psi(u))=\lim_{n
                     \rightarrow \infty} \Gamma_n(\psi(u)).$$
Furthermore, as $\|\Gamma(A)
(\psi(u))\|=\|\psi(u)\|=1$, $\Gamma(A)$
extends
to an isometry between
                     $\otimes_{i=1}^{\infty}\Gamma_s(\C e_i)$ and
                     $\otimes_{i=1}^{\infty}\Gamma_s(\C f_i)$ ((i.
e.). between $\Gamma_s(\overline{G_1})$
                     and $\Gamma_s(\overline{G_2})$).

                     Also by defining, for $u
                     \in \otimes_{i=1}^k \Gamma_s(\C f_i)$, $$
                     \Gamma^{\prime}(\psi(u)) ~=~ \lim_{n \rightarrow
\infty}(V_{{\lambda_1}^{-1}}
                     \otimes V_{{\lambda_2}^{-1}} \otimes \cdots
V_{{\lambda_k}^{-1}}
\psi(u))\otimes \bigotimes_{j=k+1}^n
V_{{\lambda_j}^{-1}} \Phi_1 \otimes
                     \Phi_1
                     \otimes \Phi_1 \otimes \cdots, $$
and by
using same arguments, we may conclude
                     that $\Gamma^{\prime}$ extends
                     to
                     an isometry between $\Gamma_s(G_2)$ and
$\Gamma_s(G_1)$, and that
                     $$\Gamma^{\prime}~=~
{\Gamma(A)}^*~=~\Gamma({A}^{-1}).$$ Hence
                     $\Gamma(A)$ is a unitary operator.

                     Clearly, as we may
conclude from equations \ref{weyl} and
\ref{weyl1}, the relations \ref{SHL1}
and \ref{SHL2} are satisfied.

Now to prove the
uniqueness, suppose there exists
                     another unitary operator $\Gamma^{\prime}$
satisfying \ref{SHL1} and
                     \ref{SHL2}, then $\Gamma^{\prime}
{\Gamma(A)}^{-1}$ commutes with all
                     Weyl operators $W(u)$. As the Weyl
representation is irreducible, we
                     conclude that $\Gamma^{\prime}=c
\Gamma(A)$, where $c$ is a complex
                     scalar of unit modulus. But the relation
\ref{SHL2}
                     implies that
$\Gamma^{\prime}=\Gamma(A)$.

\bigskip

(ii) Note that in the course of proving (i)
we have also proved
that $\Gamma(A)^{-1}=\Gamma(A^{-1})$.

Now, to prove \ref{shlrep1}, first notice, again by using the
irreducibility of the Weyl representation, that $\Gamma(AB)=c
\Gamma(A) \Gamma(B)$, for a complex number $c$ of modulus $1$.
Also it is clear from the construction that when $U$ is a unitary
operator, $\Gamma(U)=Exp(S_U)$. This is clear because $Exp(S_U)$
satisfies both
the relations \ref{SHL1} and \ref{SHL2} (note that all second quantised
operators takes the vacuum vector to the vacuum vector). It is also
easy to verify that the
relation \ref{shlrep1} is satisfied
when either $A$ or $B$ is a unitary operator (Consider equation
\ref{SHL2} and that the vaccum vector is fixed by $Exp(U)$). Hence, by
using the above fact and polar
decomposition, we may assume, without loss of generality, that
$G_1=G_2=G_3$ and that $A,B \in
\S(G)$.

We basically need to prove that $\langle \Gamma(A) \Gamma(B) \Phi,
\Phi\rangle >0$, where $\Phi$ is the vacuum vector in
$\Gamma_s(\overline{G})$.

We apeal to Proposition 22.6 in \cite{KRP} (page 166)
for the validity of the relation
$\Gamma(AB)=\Gamma(A)\Gamma(B)$, when $G$ is finite
dimensional. (The construction given in that proposition
and the
construction of
$\Gamma$ in Part (i) of this proposition are same as they both
satisfy the relations \ref{SHL1} and \ref{SHL2}.)

Let $\{A_n\}$ ( resptly. $\{B_n\}$)  $\subset \S(G)$ be
a sequence of operators,
such
that $I-A_n$ (resptly. $I-B_n$) is a finite rank operator for each
$n$, approximating $I-A$ (resptly. $I-B$), and that
$\Gamma(A_n)$
(resptly. $\Gamma(B_n)$ ) converges
strongly to $\Gamma(A)$ (resptly. $\Gamma(B)$). It is clear that
such a sequence exists from the construction. Note that $A_n$
(and similarly $B_n$ also) is a direct sum of an invertible
positive operator on the $Range(I-A_n)$($= Ker^{\perp}(I-A_n)$)
and the identity operator on $Ker(I-A_n)$, for each $n$. Let us
define for each $n$, $G_n=Span[Range(I-A_n), Range(I-B_n)]$, then
$G_n$s are finite dimensional subspaces of $G$. The following
relations
\begin{eqnarray*}
I-A_nB_n & =& I-A_n -A_n(B_n-I)\\
 &= & I-B_n -(A_n-I)B_n
\end{eqnarray*}
imply that $I-A_nB_n$ is a finite rank operator, and that
$Range(I-A_nB_n)\subset G_n$, $Ker^{\perp}(I-A_nB_n)\subset G_n$,
for each $n$. Hence we have
$A_n=P_nA_nP_n \oplus P_{[n}$,
$B_n=P_nB_nP_n \oplus P_{[n}$,
$A_nB_n=P_nA_nB_nP_n \oplus P_{[n}$, where $P_n$
and $P_{[n}$ are the
projections onto $G_n$ and $G_n^{\perp}$ respectively.
Also it is clear from the
construction of $\Gamma$ that we also have
$\Gamma(A_n)=\Gamma(P_nA_nP_n) \otimes I_{[n}$,
$\Gamma(B_n)=\Gamma(P_nB_nP_n) \otimes I_{[n}$,
$\Gamma(A_nB_n)=\Gamma(P_nA_nB_nP_n) \otimes I_{[n}$, where
$I_{[n}$ is the identity operator on
$\Gamma_s(\overline{G_n^{\perp}})$.
Therefore, as $G_n$ is finite dimensional, we may conclude
that
$\Gamma(A_nB_n)=\Gamma(A_n) \Gamma(B_n)$,  and hence that
$\langle \Gamma(A_n)\Gamma(B_n)\Phi, \Phi\rangle >0$, for each $n$.
The strong convergence of both $\Gamma(A_n)$ and $\Gamma(B_n)$ implies
that $\Gamma(A_n)\Gamma(B_n)$ converges weakly to $\Gamma(A)\Gamma(B)$.
Now
it
follows that $\langle \Gamma(A)\Gamma(B) \Phi, \Phi\rangle  >0$,
and the proof of part (ii) of the proposition is complete.

\bigskip

(iii) Suppose let $\{T_n\} \subset \S(G, G)$ converges strongly to $T
\in S(G, G)$ and $(T_n^*)^{-1}$ converges strongly to $(T^*)^{-1}$. First
we
note that the bounded set
$\overline{\{\Gamma(T_n)\}}$ (the closure is taken with respect to the
weak topology) is compact with respect to the weak
topology(Weak operator topology and weak$^*$ topology coincide on
bounded sets). Also we know
any compact $T_2$ space is metrizable, and hence the above set is
sequentially compact. So we get a convergent subsequence
$\Gamma(T_{n_k})$, say converging weakly to $V \in
B(\Gamma_s(\overline{G}))$. To prove $\Gamma(T_n)$ converges weakly to
$\Gamma(T)$, it is enough if we prove that $V = \Gamma(T)$.(This would
mean that
every subsequence of $\Gamma(T_n)$ has a further subsequence, which
converges weakly to $\Gamma(T)$, which means $\Gamma(T_n)$ converges
weakly to $\Gamma(T)$.

First we conclude, from the strong continuity of the Weyl representation
that $W(S_{T_n}x)$ converges strongly to $W(S_{T}x)$ for all $x \in
\overline{G}$. This basically means
that $\Gamma(T_n) W(x) \Gamma(T_n)^*$ converges strongly to $\Gamma(T)
W(x) \Gamma(T)^*$ for all $x \in
\overline{G}$.

We have that $$\langle \xi, \Gamma(T_{n_k})^*\Gamma(T_{n_k}) W(x)
\Gamma(T_{n_k})^*
\eta\rangle= \langle \Gamma(T_{n_k}) \xi, \Gamma(T_{n_k}) W(x)
\Gamma(T_{n_k})^*\eta \rangle$$ converging to $\langle V, \Gamma(T) W(x)
\Gamma(T)^* \eta \rangle$. We also have that $W(x) \Gamma(T_{n_k})^*$
converges weakly to $W(x)V^*$, and so we conclude that $$W(x)V^*=
V^*\Gamma(T)W(x) \Gamma(T) ~~\forall x \in \overline{G},$$ which implies
that $V^*\Gamma(T)$ commutes with all operators in $\B(\overline{G})$.
We conclude that $V$ is a scalar multiple of  $\Gamma(T)$. But by
the
fact that $V$ is the weaklimit of $\Gamma(T_{n_k})$, it follows that
$$\langle
V \Phi, \Phi\rangle > 0.$$ Hence we conclude that $V = \Gamma(T)$, and the
proof of the
theorem is
over.
\qed

\begin{Remark}
The  generalised version of Shale's theorem as presented here
for two real Hilbert spaces (part (i)) can also be proved  using the original
Shale's theorem and 
polar decomposition. That is if $A =
U(A^*A)^{\frac{1}{2}}$, then we can define, $$\Gamma(A) =
Exp(S_U)\Gamma((A^*A)^{\frac{1}{2}}),$$ where
$\Gamma((A^*A)^{\frac{1}{2}})$ is defined by original 
Shale's theorem. But  we required the 
details of the construction of $\Gamma $ in proving part (ii) of the Theorem. 
\end{Remark}

\begin{Remark}
It is clear from the construction (also a fact we have used in proof of
(ii) in Theorem \ref{shales} ), that if $A \in \S(G_1, G_2)$ and $A^\prime
\in \S(G_1^\prime,
G_2^\prime)$, then $A \oplus A^\prime \in \S(G_1 \oplus G_1^\prime,
G_2\oplus G_2^\prime)$ and $$E_2 (\Gamma(A) \otimes \Gamma(B))
E_1^*=\Gamma(A\oplus B),$$ where $E_i$ is the canonical unitary operator
between $\Gamma_s(G_i \oplus G_i^\prime)$ and $\Gamma_s(G_i) \otimes
\Gamma_s(G_i^\prime)$, $E_i(e(x\oplus y))=e(x) \otimes e(y)$, for $i=1,2$.
\end{Remark}

Our aim is to get a product system out of
what is called a `sum system'. First we define a notion of sum system
using a one parameter family of real Hilbert spaces. Later as a particular
case we will define sum system as a two parameter family  of Hilbert
spaces,
and consider only that definition throughout this paper. This 
definition is
analogous to the definition of a product system, where the tensors are
replaced by
directsums, and unitaries by our special invertible operators, which are
Hilbert-Schmidt perturbation of a unitary operators.

\begin{Definition}\label{1sst}
A sum system is a one parameter family of real Hilbert spaces $\{G_t\}_{t
\in (0,\infty)}$, together with operators $B_{s,t} \in \S(G_s \oplus G_t,
G_{s+t})$ satisfying the following axioms of associativity and
measurabilty.

\medskip
\noindent(i) (Associativity) For any $s_1, s_2, s_3 \in (0,\infty)$
$$B_{s_1, s_2 + s_3}( 1_{G_{s_1}} \oplus B_{s_2 ,
s_3})=  B_{s_1+ s_2 , s_3}( B_{s_1 ,
s_2} \oplus 1_{G_{s_3}}).$$

\medskip
\noindent (ii) (Measurability) There exists a countable set $G^0$ of
sections $t\in
R \rightarrow x_t \in G_t$ such that $ t  \mapsto \langle x_t^n,
x_t^m\rangle$ is
measurable for any two $x^n, x^m \in G^0$, and the set $\{x_t^n: x^n \in
G^0\}$ is a total set in $G_t$ for all $ t \in (0,\infty)$.
Further it is also assumed that the following maps $$t\in \R  \mapsto
B_{t,1-t}(x^m_{t}
\oplus 0) \in G_1, ~ t \mapsto B_{t,1-t}(0
\oplus x^n_{1-t})$$ are measurable for any fixed $n,m \in \N$.
\end{Definition}

Given a Sum system $(G_{t}, B_{s,t})$, define $$H_t =
\Gamma_s(\overline{G_t}),~ U_{s,t} = \Gamma(B_{s,t}),$$ where the Hilbert
spaces $\Gamma_s(\overline{G_s}) \otimes \Gamma_s(\overline{G_t})$ and
$\Gamma_s(\overline{G_s} \oplus \overline{G_t})$ are
identified,  using the canonical unitary operator taking $e(x) \otimes
e(y)$ to $e(x \oplus y).$ Now we have produced a product system from a
given
sum system.

\begin{Theorem}\label{pdct1}
$(H_t, U_{s,t})$, defined as above, is a product system.
\end{Theorem}

{\em Proof:}
The associativity property follows from the associativity of the sum
system, and from statement (ii) of theorem \ref{shales}. So we basically
have to prove the axiom of measurability.

To prove the measurability axiom, we use the group of
unitary operators $\{\tau_t\}$ on $H_1$, defined in \cite{Vol}. Let
$\pi_t$ be
the unitary map between $H_{1-t} \otimes H_t$ and $H_t \otimes
H_{1-t}$ given by $\pi_t(x_{1-t} \otimes x_t)= x_t \otimes
x_{1-t}.$ Then define for each $t \in (0,1)$, a unitary
operator on $H_1$,  by $\tau_t=U_{t, 1-t}
\pi_t U_{1-t, t}^*,$ and we set $\tau_1 = 1_{H_1}$, and
$\tau_{t+k}= \tau_t$ for  any $ k \in \Z$.

It is proved in \cite{Vol} that $\{\tau_t\}_{t \in \R}$
forms an one parameter unitary group (see Proposition 2 in
\cite{Vol}).  It is also proved in \cite{Vol} that all measurable
structures on a given algebraic product system leads to isomorphic
product systems, and an algebraic product system admits a measurable
structure if and only if the unitary group $\{\tau_t\}$ is
continuous (theorem 51 in \cite{Vol}). Therefore we prove that
$\{\tau_t\}$ is strongly continuous.

Define a group of operators, $T_t$ on the real Hilbert space $G_1$, by
$T_t = B_{t, 1-t}
\sigma_t B_{1-t, t}^{-1}$, where $\sigma_t: G_{1-t} \oplus G_{t} \mapsto
G_{t} \oplus G_{1-t}$, is the unitary operator, defined by $\sigma_t(x
\oplus y)= y \oplus x$, for $t \in (0,1)$. Also set $\T_1 = 1_{G_1}$, and
$T_{t+k}= T_t$ for  any $ k \in \Z$. The fact that $T_t$ is a group
can be checked in same way for $\tau_t$.  Now it is easy to check that
$T_t \in \S(G_1, G_1)$, and using statement (ii) of theorem \ref{shales}
it
is also clear that $\tau_t = \Gamma(T_t)$. As adjoint of a 
strongly continuous semigroup is again a strongly continuous 
semigroup (see Theorem 4.3 of \cite{Gold}), suppose if we 
prove that the group $\{T_t\}$ is strongly continous, then the group 
$\{(T_t^*)^{-1}\}$ is also strongly continuous. Then this would imply 
the weak continuity, hence the strong continuity, of the unitary 
group
$\{\tau_t\}$, by the statement (iii) in theorem \ref{shales}. So we prove
the strong continuity of $\{T_t\}$. This is
equivalent to prove the strong measurability of ${T_t}$ (see \cite{HP},
part two, chapter X), and by the definition of $T_t$ it is enough to
prove the measurability for $t \in (0,1)$.

Let us assume that the set of all measurable sections is indexed by
$\N$. Define $y_t^k \in G_1$, for $k \in\N, t \in (0,1)$ by
$$ y^{2k-1}_t = B_{1-t, t}(x^k_{1-t} \oplus 0), ~y^{2k}_t = B_{1-t, t}(
0 \oplus x^k_t).$$ Then the invertibility of $B_{1-t, t}$ implies that
the set $\{y_t^k\}_{k \in \N}$ is a linearly independent and total set in
$G_1$.  Let $\xi_k^t$ be the Gram-Schmidt orthogonalisation of $y_k^t$, i.e.
$$\xi^1_t = \frac{y^k_t}{\|y_k^t\|}, ~ {}^\prime\xi^{k+1}_t = y_t^{k+1} -
\sum^k_{i=1} \langle y^{k+1}_t, \xi_t^i\rangle \xi_t^i, ~\xi^{k+1}_t =
\frac{{}^\prime\xi^{k+1}_t}{\|{}^\prime\xi^{k+1}_t\|}, ~\mbox{for} ~ k
\in \N.$$  The measurability axiom of the sum system says that map $t
\mapsto y^k_t$ is measurable for $ k \in \N$. It is an easy verification,
using induction, to see that the map $t \mapsto \xi^k_t$ is also
measurable.

We need to prove that the map $ t \mapsto T_t x^n_1$ is measurable, for
any fixed $n \in \N$. Now, $$T_t(x_1^t ) = \sum_k \langle x_1^n,
\xi_t^k\rangle T_t\xi^k_t.$$ We basically need to prove that the map
$t\mapsto T_t\xi^k_t$ is measurable. Notice  that $T_t B_{1-t, t} =
B_{t, 1-t}$. Using the fact that $ t \mapsto \xi_t^k$ is measurable
and induction, we may conclude that the map $ t \mapsto T_t \xi^k_t$ is
measurable.

We have proved the measurability axiom for the product system, and that
$(H_t, U_{s,t})$ forms a product system.
\qed

\bigskip

We call $(H_t, U_{s, t})$ as the exponential of the sum system $(G_t, B_{s,t})$
or as the product system arising out of this sum system.
Next we define the notion of isomorphism for sum systems.

\begin{Definition}
Two sum systems $(G_t, B_{s,t})$ and $(G^\prime_t, B^\prime_{s,t})$
are
said to be isomorphic if there exists an operator $A_t \in \S(G_t,
G_t^\prime)$ for each
$t \in (0,\infty)$, satisfying $A_{s+t}B_{s,t}= B_{s,t}^\prime (A_s \oplus
A_t).$
\end{Definition}

Clearly, by statement (ii) in theorem \ref{shales}, if two sum systems are
isomorphic, then the corresponding product systems are also isomorphic, where
the isomorphism between the product systems are implemented by $\Gamma(A_t),~t
\in (0,\infty)$. It is not clear as to whether the converse is true. 
 Next we define a sum system given by a two parameter family of
Hilbert spaces, and a semigroup of shift operators. Two parameter systems are
more convenient. All our examples will be of this kind.

\begin{Definition}\label{sumsystem}
A two parameter sum system is a two parameter family of real Hilbert spaces
$\{G_{(s,t)}\}$ for $0 <s < t< \infty$ all embedded into a
single linear space $G^0_{(0, \infty)}$,
satisfying $G_{(s,t)} \subset
G_{(s^{\prime}, t^{\prime})}$ if the interval
$(s,t)$ is contained in the
interval $(s^{\prime}, t^{\prime})$, together
with a one parameter
semigroup $\{S_t\}$, of linear maps on
$G^0_{(0, \infty)}$ for $t \in (0,
\infty)$ such that

\medskip
\noindent
(i) $S_s|_{G_{(0,t)}} \in \S(G_{(0,t)} , G_{(s, s+t)})$.

\medskip
\noindent
(ii) $G_{(0, s+t)}= G_{(0,s)} \uplus G_{(s, s+t)}$ for all
$s,t$.

\medskip
\noindent
(iii) The semigroup $\{S_t\}$ is `locally' strongly continuous, (i.e.)
for any $x
\in G_{(a,b)}$, $a,b \in (0,\infty)$, $S_tx$ converges to $x$, as
$t \rightarrow 0$, where the convergence takes place in a bigger Hilbert
space, $G_{(a, b+\epsilon)}$, for some $ \epsilon >0$.
\end{Definition}

Notice that the condition (iii) in the above definition actually
implies that $S_t$ converges strongly to $S_{t_0}$, if $t \rightarrow
t_0$, due to the semigroup property. We may assume that
$G^0_{(0,\infty)}=\cup_{t >0}G_{(0,t)},$ and
define $G_{(0,\infty)}= \overline{G^0_{(0,\infty)}}$, as the Hilbert space
completion. (The problem is we may not be able to extend the semigroup $S_t$ to
$G_{(0,\infty)}$.)


Let $(G_t, B_{s,t})$ be a (one parameter) sum system such that $B_{s,t}|_{G_s}$
is an isometry, for $s,t \in (0,\infty)$. Then $(G_t, B_{s,t})$ can
be shown
to be isomorphic to a sum system given by a two parameter family, in the
following way.
We can define $G_{(0,\infty)}$ as the inductive limit of
the Hilbert spaces $G_s$. That is define $$\tilde{G}_{(0,\infty)} ~=
~\bigcup_{t> 0}G_t.$$
Define an equivalance relation on $\tilde{G}_{(0,\infty)}$ by the
following,
for $x \in G_s$ and $y \in G_t$ and $t >s$,  $x
\sim y ~~\mbox{if}~ B_{s,t-s}x =y.$ The associativity axiom implies that (by
taking $s_1=s,~ s_2= t-s,~ s_3= t^\prime -t$) $$(B_{t, t^\prime -t}B_{s,
t-s})|_{G_s} = B_{s, t^\prime-s}|_{G_s}.$$ Hence if $B_{s, t-s}x= y$ and
$B_{t,t^\prime-t}y =z$, then $B_{s, t^\prime -s}x =z$. So we have an
equivalence relation.
Define $$G^0_{(0,\infty)} ~=~ \tilde{G}_{(0,\infty)}/\sim,~ i_t:G_s \mapsto
G^0_{(0,\infty)},~ i_t(x) = [x]. $$ We can define $$\lambda [x] = [\lambda x],~
[x] + [y] = [x+y],$$ where the sum $x+y$ is taken by embedding $x$ and $y$ in
a common bigger Hilbert space(which will be again consistent by the
associativity axiom). If we define $\|[x]\|=\|x\|$ (which is well
defined due to the isometric assumption on $B_{s,t}|_{G_s}$), then $i_t$
is an embedding of
$G_t$ into $G^0_{(0,\infty)}$. Define
$$G_{(0,\infty)}= \overline{G^0_{(0,\infty)}}, ~ G_{(0,t)}=i_t(G_t),
~\mbox{and clearly}~ G_{(0,\infty)}=\overline{\bigcup_{t
>0}G_{(0,t)}}.$$  For $ [x] \in G_{0,t} \subset G^0_{(0,\infty)}$ define
$$S_s([x]) =
[B_{s,t}x],~\mbox{for} ~ s,t \in (0,\infty).$$
It can be checked, again by using the associativity axiom, that the map
$S_s$ is well defined and that $\{S_t\}$ forms a semigroup also. Finally
the the strong continuity of $\{S_t\}$ will follow from the measurability
axiom.

\begin{Remark}
It is not clear as to whether a general (one parameter) 
sum system is  isomorphic to a sum system such that 
$B_{s,t}|_{G_s}$ is isometric for all $ s \in (0,\infty)$. 
\end{Remark}


Given a two parameter sum system $(G_{(s,t)}, S_t)$ we get a
one parameter
sum system by defining, $$G_t = G_{(0,t)}, ~ B_{s,t}(x_s \oplus y_t) =x_s
+ S_s y_t.$$ Then clearly the associativity axiom is satisfied by $(G_t,
B_{s,t})$, due to the semigroup property of $S_t$. The measurability axiom may
be proved as follows.

Let $P_t$ denote the orthogonal projection from $G_{(0,\infty)}$ onto
$G_{(0,t)}$, and let $\{x^n\}_{n \in \N}$ be any orthonormal basis for
$G_{(0,\infty)}$. Define
$x^n_t= P_tx^n$. Then clearly $x^n_t$ is a countable total set in
$G_{(0,t)}$,
for each $t \in (0,\infty)$. Also clearly $P_t \uparrow I$ as $t \rightarrow
\infty$, and $P_t \uparrow P_{t_0}$ as $ t \uparrow t_0$ for any $t_0 \in
(0,\infty)$. Hence the map $t\mapsto P_t x_n$ is measurable, and in particular
the map $$t \mapsto \langle x_t^n, x_t^m\rangle =\langle P_t x^n, x^m\rangle$$
is measurable for any $n, m \in \N$. Again clearly $$ t \mapsto B_{t,
1-t}(x_t^n\oplus 0) = x_t^n$$ is measurable. So we only have to prove
that the
map $$ t \mapsto B_{t, 1-t}(0 \oplus x^ m_{1-t}) = S_t(x^m_{1-t})$$ is
measurable. Denote the above map by $f(t) = S_t(x^m_{1-t})$. Now define
for $ k \in
\N$, a function $f_k:(0,\infty) \mapsto
G_{(0,1)}$ by $$f_k(t) = S_t(x_{1-\frac{l}{k+1}}) ~~~~\mbox{if}~ t \in
(\frac{l}{k}, \frac{l+1}{k}),~l=0,1, \cdots
k-1.$$
Clearly
the function $f_k $ is measurable for each $ k\in \N$
due to the strong continuity of $S_t$, and $f_k$ converges to $f$ pointwise, as
$x_{t_n} \rightarrow x_t$ if $t_n \uparrow t$. Now the measurability of the
function $f$ proved.


\begin{Remark} In this construction of a one parameter sum system out of a two
parameter sum system we have used the map: $B_{s,t}(x_s\oplus
y_t)= A_{s,t}(x_s\oplus S_sy_t)$ where $A_{s, t}$ is the map
$x\oplus y\mapsto x+y$. Instead of this $A_{s,t}$ we could have
used any map  $A_{s,t} \in \O(G_{(0,s)} \oplus G_{(s,s+t)}, G_{(0,s+t)})$
and this has no effect on the product system arising out of the
sum system because of the following Lemma.
\end{Remark}

\begin{Lemma}\label{ind}
Let $G_1, G_2 \subset G$ and $A \in \O(G_1\oplus G_2,G)$. Define
$H_i=\Gamma_s(\overline{G_i}),~i=1,2~,~
H=\Gamma_s(\overline{G})$, and a unitary operator $V$
between $H_1 \otimes H_2$ and
$H$, by
\begin{eqnarray}\label{unitry0}
V & = & \Gamma(A)~ E~
(\Gamma(A^{-1}|_{G_1}) \otimes
\Gamma(A^{-1}|_{G_2})),
\end{eqnarray}
where $E$ is the canonical
unitary operator between
$\Gamma_s(G_1) \otimes \Gamma_s(G_2)$ and $\Gamma_s(G_1 \oplus
G_2)$ and $\Gamma(A)$ provided by Shale's theorem. Then $V$ does not
depend on the particular choice of $A$.
\end{Lemma}

{\em Proof:} Let $A_1, A_2 \in  \O(G_1\oplus G_2,G)$, then clearly
$A_2^{-1}A_1 \in \O(G_1 \oplus
G_2, G_1 \oplus G_2)$, which basically means that $A_2^{-1}A_1$
splits into direct sum of
operators, i. e. $A_2^{-1}A_1=(A_2^{-1}A_1)|_{G_1} \oplus
(A_2^{-1}A_1)|_{G_2}$. Hence, we conclude that
$$\Gamma(A_2^{-1}A_1)= E\Gamma((A_2^{-1}A_1)|_{G_1}) \otimes
\Gamma((A_2^{-1}A_1)|_{G_2})E^*,$$ where $E$ is the canonical
unitary operator between $\Gamma_s(G_1) \otimes \Gamma_s(G_2)$ and
$\Gamma_s(G_1) \oplus \Gamma_s(G_2)$. By applying
relations \ref{shlrep} and \ref{shlrep1} we
may conclude that
$$\Gamma(A_2)^*\Gamma(A_1)=E\Gamma(A_2^{-1}|_{G_1})\Gamma(A_1|_{G_1})
\otimes
\Gamma (A_2^{-1}|_{G_2}) \Gamma(A_1|_{G_2})E^*,$$ which would
prove that $V$ does not depend on the particular choice
of $A$.
\qed

Here after we normally take only two parameter sum systems and
we construct the one parameter sum system, and then the product system
from it using the map $A_{s,t}(x_s\oplus y_t)= x_s+y_t.$

To begin with we present two sets   of examples for sum systems.
First one was given by Arveson producing the
type $I$ exponential product system, when the sum system comes from 
usual $L^2$ on intervals. The
other one
is the example of Tsirelson, producing type $III$ product
system where the $L^2$ spaces are completed with 
respect to a different
inner product coming from  carefully chosen positive definite
kernels. In the next
Section we will see that under some simplifying assumptions 
only type $I$ and type $III$ arise as product systems of sum systems.
 In
particular
it seems to be impossible to produce type $II$ product systems from a
sum system.

\begin{Example}
Let $G_{(a,b)}=L^2((a, b), K)=\{f :(a,b) \rightarrow K: \int\|f\|^2 < \infty\}$,
where $K$ is a separable Hilbert space, and
$S_t$ be the usual shift $S_t(f)(s)= f(s-t)$. Then exponential of this sum system is
 the
exponential or Fock product
system of Arveson, given in \cite{Arv}. These are completely classified by the dimension of
$K$. 
\end{Example}

\begin{Example}\label{Tsirelson}
In \cite{Ts}, Tsirelson defines a scalar product on $L^2(a,b)$, given by
\begin{eqnarray}\label{sclrprdt}
\langle f,
g\rangle = \int \int f(s)g(t) B(s-t)ds dt,
\end{eqnarray}
where $B \in L^1(\R)$ is
continuous and
positive definite. Let $G_{(a,b)}$  be the completion of $L^2(a,b)$ with respect to
this inner product and
let $S_t$ be the usual shift $S_t(f)(s)= f(s-t)$ extended. 
Then $\{ S_t\}$ is a strongly continuous
semigroup of isometries. It
is also assumed that $B$ satisfies the following property
\begin{eqnarray}\label{condition}
\exists ~\epsilon >
0~\mbox{such that}~
\forall t \in (0, \epsilon)~B(t)= \frac{1}{t \ln^{\alpha} (\frac{1}{t})}
\end{eqnarray}
and
the function $B$ is positive, decreasing and convex. With this assumption it is
proved that the map $x \oplus y \rightarrow x+y$ is in $\S(G_{(0,s)}
\oplus G_{(s,
s+t)}, G_{(0, s+t )})$ (see proposition 9.9 in page 48, \cite{Ts}). So $(G_{(a,b)},
S_t)$ forms a sum system. We will prove in the
next two sections, that the corresponding product systems (for different $\alpha$
in  the condition \ref{condition}) are unitless and non-isomorphic.
\end{Example}

Before ending this section, we prove some facts regarding sum systems.
First we prove that single points does not matter in a sum system, in the
following sense. Given a sum system, we can naturally associate a real
Hilbert space to any given interval. It does not matter whether the end
points of the interval are included or not. This basically follows from
the our assumption that the shift semigroup is strongly continuous.

We prove two easy lemmas before that. The first one is about the uniform
boundedness of the shift semigroup over any finite interval, which is
a well known fact for any strongly continuous semigroup on
Banach spaces.

\begin{Lemma}
For any $a, b, s_1, s_2 \in (0, \infty)$, $$\sup_{t \in (s_1,
s_2)}\|S_t|_{G_{(a,b)}}\| < \infty.$$
\end{Lemma}

{\em Proof:} We use the uniform boundedness principle.(We may consider
the family $\{S_t\}_{t \in (s_1, s_2)}$, as operators between
the two Banach spaces, $G_{(a,b)}$ and $G_{(a+ s_1,b+ s_2)}$.)   So we
only need to prove that for any $x \in G_{(a,b)}$, $$\sup_{t \in (s_1,
s_2)}\|S_t|_{G_{(a,b)}}x\| < \infty .$$  Suppose there
exist a sequence $\{t_n\} \subset (s_1, s_2)$ converging to $t \in
(s_1,
s_2)$, and $\|S_{t_n}x\| \geq n$, for each $n \in \N$. But then
$S_{t_n}x$ can not converge to
$S_tx$, which contradicts the strong continuity assumption of $\{S_t\}$.
\qed

\bigskip

 From here onwards we denote the restriction of the shift semigroup as
just $S_t$, unless there is any confusion.

\begin{Lemma}
For any $ x \in G_{(0,1)}$, $S_tT_tx$ converges to $x$ and $S_1x$, as
$t$ tends to $0$ and $1$ respectively, where $T_t$
is the
semigroup which is already defined by $$T_t(x) = x_{1-t}^\prime + S_{1-t}
x_t,~~\mbox{if}~~x = x_t + S_tx_{1-t}^\prime, ~\mbox{for}~ x_t \in
G_{(0,t)},~x_{1-t}^\prime \in G_{(0, 1-t)}.$$
\end{Lemma}

{\em Proof:}
We have, $$\|S_tT_t(x)
-x\|\leq \|S_t\|\|T_tx -x\| + \|S_tx
-x\|.$$
Similarly
we also have, $$\|S_tT_t(x)
-S_1x\|\leq \|S_t\|\|T_tx -x\| + \|S_tx
-S_1x\|.$$
\qed

\begin{Proposition} For $t \in (0,1)$, let $ x \in G_{(0,1)}$ be
such that $x = x_t +
S_tx_{1-t}^\prime$ for $ x_t \in
G_{(0,t)},~x_{1-t}^\prime \in G_{(0, 1-t)}.$ Then $x_t$ and
$S_tx_{1-t}^\prime$
converges to $0$, as $t$ tends to $0$ and $1$ respectively.
\end{Proposition}

{\em Proof:} We have $$x = x_t + S_t(x_{1-t}^\prime)=x_t + S_tT_t(x)
-S_1x_t.$$ Hence by the above lemma $(I -S_1)x_t$ converges to $0$ as $t
\rightarrow 0$. Similarly we also have $$S_tT_t(x)=S_t(x_{1-t}^\prime +
S_{1-t}x_t)= (S_tx_{1-t}^\prime + S_1x - S_1(S_tx_{1-t}^\prime).$$ Again
the above lemma implies that $(I-S_1)S_tx_{1-t}^\prime$ converges to $0$.

The map $(I-S_1) : G_{(0,1)} \mapsto G_{(0,2)}$, is clearly injective,
and hence a bijection between $G_{(0,1)}$ and its range. Notice that the
proof of the Proposition is over if we prove that the inverse is bounded.
To  prove that first notice that the map between $G_{(0,1)} \mapsto
G_{(0,1)} \oplus G_{(0,1)}$ given by $x \mapsto x \oplus -x$, is a
bijection between $G_{(0,1)}$ and its range, with a bounded inverse. The
remaining part of the
proof follows from the property (ii) in the definition of a sum system.
\qed

\begin{Corollary}\label{Gt+}
Let $(G_{(a,b)}, S_t)$ be a sum sytem, then $$G_{t+}=\bigcap_{s>0} G_{(t,
t+s)}= \{0\},~G_{t-}=\bigcap_{s>0} G_{(s, t)}= \{0\}.$$
\end{Corollary}

{\em Proof:} As each $S_t$ in the shift semigroup is a bijective map, it
is enough if
we prove that $$G_{0+}=\{0\}=G_{1-}.$$ Suppose $x \in G_{0+}$, then the
decomposition in the above proposition becomes $x_t =x$ and
$x_{1-t}^\prime =0$ for any $ t \in (0,1)$. Hence $x_t =0$. In an exactly
similar way, from the other part of the above proposition, we may
conclude that $G_{1-}=\{0\}$.\qed

\section{Invariants}\label{non-isomorphic}

In this section we get an invariant for any product system
constructed out of a sum system. The invariant we get is same
as the one got by Tsirelson in \cite {Ts}, but we prove it in
our setup. Also the proof turns out to be more direct and simple.

Let $(H_t, U_{s,t})$
be any product system. Associate for any closed interval
$[s,t] \subset [0,1]$, a von
Neumann algebra defined by $${\cal A}_{[s,t]}~=~ U_{s,t,
1}~(1_{H_s} \otimes
\B(H_{t-s})\otimes
1_{H_{1-t}}) U_{s,t,1}^*,$$ where $U_{s,t,1}$ is the canonical unitary
operator
between the Hilbert spaces $H_s
\otimes H_{t-s} \otimes H_{1-t}$ and $H_1$, determined uniquely by the
associativity
of the
product system.
We define an elementary set to be a subset of $[0,1]$, which is
disjoint union of finite number of closed intervals. We denote by
${\cal F}^e={\cal F}^e_{[0,1]}$ the collection of all elementary sets in
$[0,1]$.
For an
elementary set $E= \sqcup_{i=1}^n[s_i,t_i]$, define the associated
von Neumann  algbra to be the von Neumann algebra generated by all the
von
Neumann algebras associated with the individual intervals, i.e.
$${\cal A}_E~=~\bigvee^{n}_{i=1}{\cal A}_{[s_i,t_i]}.$$

We define the concept of $\liminf$ for a sequence of von Neumann algebras
as follows.

\begin{Definition} For a sequence of von Neumann algebras ${\cal A}_n$ we
define $\liminf {\cal A}_n$ as the von Neumann algebra generated by limits
of all subsequences $\{T_{n_k}\}$, of any sequence $\{T_n\}$ such that
$T_n \in {\cal A}_n$, where the
limit is
taken in the weak operator topology.
\end{Definition}

Clearly the set of all sequences of elementary sets $E_n$ such that
$\liminf
{\cal A}_{E_n} = \C$, is an invariant of the product system under
isomorphisms. From this observation we get the invariants for the product
system, given in terms of the sum system, by Tsirelson.

When the product system arises from a sum system, we define
$G_E$ for $E \in {\cal
F}^e$, to be the Hilbert space
generated by all Hilbert spaces corresponding to the individual
intervals. We will talk about $G_E$ and ${\cal A}_E$, only when $E$ is an
elementary
set, so it does not matter whether the intervals are closed or not, due
to Corollary  \ref{Gt+} in the previous section.

In order to get the invariants for the product systems, arising from a
sum system, we also make some definitions of $\liminf$ and $\limsup$ of
subspaces of a Hilbert space. We will be making use of these concepts in
the next section also. The definitions are
same as in \cite{Ts}.

\begin{Definition} Let $G$ be a real Hilbert space, and
$\{G_n\}_{n \in \N}$
be a sequence of subspaces of $G$, then $ \liminf G_n = \{x \in G: x=
\lim x_n, ~ x_n \in G_n\}.$ Also we define the
$\limsup G_n$  to be the closed subspace generated by weak limits of all
subsequences of $x_n$, such that $x_n \in G_n$, , i. e.
$$\limsup G_n =\overline{span}\{x: w-\lim x_{k_n}=x,~ x_{k_n}\in
G_{k_n}, ~ \{k_n\} \subset \N \}.$$
\end{Definition}

\begin{Lemma}\label{infsup}
Let $G$ be a real Hilbert space. For any sequence of
subspaces $G_n$, $\limsup G_n = (\liminf G_n^{\perp})^{\perp}.$
\end{Lemma}

{\em Proof:}
First we will prove the inclusion $\limsup G_n \subset (\liminf
G_n^{\perp})^{\perp}.$ That is we need to prove that
$$
\liminf
G_n^{\perp}
\subset
(\limsup
G_n)^{\perp}.$$

Let
$y \in \liminf G_n^{\perp}$, that is there
exists
a sequence $\{y_n\}$ such that $y_n \in G_n^{\perp}$ and $y_n$ converges
to $y$. Also let
$x \in G$, be the weak limit of of some sequence $x_{k_n}$,
where
$x_{k_n} \in G_{k_n}$.
Then it is easy to verify that $\langle x, y \rangle = \lim \langle
x_{k_n},
y_{k_n}\rangle=0.$ This proves the required inclusion.

To prove the other inclusion, it is enough if we prove that
$$(\limsup G_n)^{\perp} \subset \liminf G_n^{\perp}.$$  Let $ y \in
(\limsup G_n)^{\perp}$, and let $y_n= P_n y \in G_n^\perp$, where $P_n$ is
the
orthogonal projection onto $G_n^{\perp}$. It is enough to prove
that $y_n$
converges to $y$, that is $y - y_n \in G_n$ converges to $0$. Note that
$$\|y-y_n\|^2=\langle y -y_n, y-y_n\rangle =\langle y, y-y_n\rangle.$$
Hence it is enough to prove that every subsequence of $y-y_n$ has a
further weakly convergent subsequence which converges weakly to $0$. As $y
-y_n$ is
bounded, every subsequence has a weakly convergent subsequence. Now
suppose
$\{y - y_{k_n}\}$ be a convergent subsequence of $y-y_n$
converging to $x$, then by definition $x \in \limsup G_n$, and hence by
our assumption $\langle y, x \rangle =0$. Now $\|y -y_{k_n}\|^2=\langle
y, y -y_{k_n}\rangle$ converges to $\langle y, x\rangle =0$.
The proof of the lemma is
over.
\qed

\bigskip

In our setup (i.e. when the product system is constructed from a
sum system), for
a set $E \in {\cal F}^e$, and $E=\bigsqcup_{i=1}^n [s_i,t_i]$, we have
$${\cal
A}_E~=~\Gamma(A_E)
(\otimes_{i=1}^n\B(\Gamma_s(\overline{G_{(s_i,
t_i)}})) \otimes_{i=0}^n
1_{\Gamma_s(\overline{G_{(t_i,s_{i+1})}})}){\Gamma(A_E)}^{-1},$$
where
we assume $t_0 =0$ and $s_{n+1}=1$, and $A_E \in \S(\oplus_{i=1}^nG_{(s_i,
t_i)}\oplus _{i=0}^n G_{(t_i, s_{i+1})}, G_{(0,1)})$ is
the map
taking $\oplus
x_{s_i,t_i}\oplus x_{t_i, s_{i+1}}$ to $\sum x_{s_i,t_i}+ \sum x_{t_i,
s_{i+1}}$.

Noting that
the von Neumann algebra $\otimes_{i=1}^n\B(\Gamma_s(\overline{G_{(s_i,
t_i)}}))$ is generated by the set of Weyl operators $\{W(x +iy):x, y \in
\oplus_{i=1}^n G_{(s_i, t_i)}\}$, it is easily seen that ${\cal
A}_E$ is
generated by the set of Weyl operators
$$
\{W(A_Ex+ i(A_{E}^*)^{-1}y): x, y \in
\oplus_{i=1}^n G_{(s_i, t_i)}\}.$$ Then, as $A_E$
(resptly. $(A^*_E)^{-1}$) is a
bijection between $\oplus_{i=1}^n G_{(s_i, t_i)}$ and
$G_E\subset G_{(0,1)}$
(resptly. $G_{E^c}^{\perp} \subset G_{(0,1)}$),
we conclude that $${\cal
A}_E~=~VNalg\{W(x +iy): x\in
G_E, ~ y \in
G_{E^c}^{\perp}\},$$ for $E \in {\cal F}^e$.

Also, using the
fact that $W(x)$ (resptly. $W(ix)$) commutes with $W(iy)$ (resptly.
with $W(y))$ when $x$ and $y$ are orthogonal vectors, and by looking at
the generators, it is easy to check that, for $E
\in {\cal F}^e$, we have
$${\cal A}_E^\prime~=~{\cal
A}_{E^c}~=~VNalg\{W(x +iy): x\in
G_{E^c}, ~ y \in
G_{E}^{\perp}\} ~~~~~~~~~~~~~~~~~~~~~~~~~(*)$$

We prove a lemma which will be used in the main theorem of this
section.

\begin{Lemma}\label{inclcommt}
Let $\{F_n\}$ be any sequence of elementary sets, then

\medskip \noindent
(i) $VNalg\{W(x+iy): x \in \limsup G_{F_n}, ~ y \in \limsup
G_{F_n^c}^{\perp}\} \subset \liminf {\cal A}_{F_n}.$

\medskip \noindent
(ii) $VNalg\{W(x+iy): x \in
\liminf G_{F_n^c}, ~ y \in \liminf G_{F_n}^{\perp}\} \subset (\liminf{\cal
A}_{F_n})^\prime.$
\end{Lemma}

{\em Proof:}
(i) Let $x \in \limsup G_{F_n}$, that
is $x=w-\lim x_{n_k}$, for some subsequence $ x_{n_k}$ such that
$x_{n_k} \in G_{F_{n_k}}.$ Then $e^{\|x_{n_k}\|^2} W(x_{n_k}) \in {\cal
A}_{F_{n_k}},$ and it is an easy verification to check that $\langle
e^{\|x_{n_k}\|^2} W(x_{n_k}) e(y), e(z)\rangle$ converges to $\langle
e^{\|x\|^2} W(x) e(y), e(z)\rangle$, for all $y,z, \in G_{(0,1)}$. Hence
we
conclude that $ W(x) \in \liminf {\cal A}_{F_n}$. Using the same argument
we may
conclude that $W(iy) \in \liminf {\cal A}_{F_n}$, for $ y \in  \limsup
G_{F_n^c}^{\perp}.$

\medskip \noindent
(ii) Let $x \in \liminf G_{F_n^c}$, that is $ x = \lim x_n$,
where $x_n \in G_{F_n^c}$. Also let $a \in \liminf {\cal A}_{F_n}$, that is
there exists
a sequence $a_{n_k} \in {\cal A}_{F_{n_k}}$, such that $a_{n_k}$ converges
in the weak operator topology to $a$. We want to prove that $W(x)$
commutes with $a$. We have that $W(x_{n_k})$ (and its adjoint
$(W(-x_{n_k})$) converges strongly to $W(x)$ (respectively to its adjoint
$W(-x)$), and that $a_{n_k}$ (and its adjoint $a^*_{n_k}$)
converges weakly to $a$ (respectively to $a^*$). Using the observation
$(*)$ above, we note that
that
$W(x_{n_k})$ and $a_{n_k}$ commutes with each other.
For any $ \xi, \eta \in H_1,$ $$\langle a W(x)\xi, \eta\rangle = \lim_k
\langle
W(x) \xi, a_{
n_k}^*\eta\rangle,$$ and $$\langle
W(x) \xi, a_{
n_k}^*\eta\rangle \leq \langle W(x_{n_k}) \xi, a_{n_k}^*\eta \rangle +
\|W(x)\xi-W(x_{n_k})\xi\|\|a_{n_k}^* \eta\|.$$ As $\|a_{n_k}^* \eta\|$ is
bounded, we get $$\langle a W(x)\xi, \eta\rangle = \lim_k \langle
W(x_{n_k}) \xi, a_{
n_k}^*\eta\rangle = \langle a_{n_k} \xi, W(-x_{n_k}) \eta\rangle.$$ Now using
the
same convergences of sequences and retracing the same arguments we may
conclude
that
$$\langle a W(x)\xi, \eta\rangle = \langle W(x) a  \xi, \eta\rangle.$$
Hence $W(x) \in (\liminf {\cal A}_{F_n})^\prime$. A similar calculation will
imply that
$W(iy) \in (\liminf{\cal A}_{F_n})^\prime$, for any $ y \in \liminf
G_{F_n}^{\perp}.$ The lemma is proved
\qed

\bigskip

The following theorem allows us to compare the invariants
through the sum system.

\begin{Theorem}\label{invariant} Let $F_n$ be any given sequence of
elementary sets,
then the following two statements are equivalent.

\medskip
\noindent
(i) $\liminf {\cal A}_{F_n} = \C.$

\medskip
\noindent
(ii) $\liminf G_{F_n^c} = G_{(0,1)}, ~ \limsup G_{F_n} = \{0\}.$
\end{Theorem}

{\em Proof:} We first prove (i)
implies (ii). We conclude using lemma \ref{infsup} and part (i) of lemma
\ref{inclcommt},
that $$VNalg\{W(x+iy): x \in \limsup G_{F_n}, ~ y \in (\liminf
G_{F_n^c})^{\perp}\} \subset \liminf {\cal A}_{F_n},$$
and clearly (i) implies (ii).

Now we prove the other implication, (ii) implies (i). Again using lemma
\ref{infsup} and part (ii) of lemma \ref{inclcommt} we have that
$$VNalg\{W(x+iy): x \in
\liminf G_{F_n^c}, ~ y \in \liminf G_{F_n}^{\perp}\} \subset (\liminf {\cal
A}_{F_n})^\prime.$$ If we assume (ii) holds, then LHS in the above inclusion
is $B(H_1)$ and the $(*)$ implies that (i) is true.
The proof of the  theorem is over
\qed

\bigskip

\begin{Remark}
The above theorem asserts that the collection of all sequence of
elementary sets
$\{E_n\}$ such that
$\liminf G_{E_n}= G_{(0,1)},$  and
$\limsup G_{E_n^c}=\{0\}$
is an invariant of the product systems.
Tsirelson has produced sequence of elementary sets satisfying
$\liminf G_{E_n}= G_{(0,1)},$  and
$\limsup G_{E_n^c}=\{0\}$,
for each
$\alpha$
but it violates the condition
$\limsup G_{E_n^c}=\{0\}$ for $\alpha^\prime \neq \alpha$.
This proves that the examples of Tsirelson are non-isomorphic for
different values of $\alpha$.
\end{Remark}

\section{Units in the product system}\label{onunits}

In this section we get a sufficient condition for the product
system, arising from
what is called as a divisible sum system, to be unitless. We
prove a necessary
condition for a unit to exist, and the sufficient
condition for the product system to be unitless is to
violate that. We first define the notion of divisibility for sum systems
and prove that
this property is satisfied by the examples of Tsirelson. All through this
section, we
assume that the restriction of the shift map $S_t|{G_{(a,b)}}$ of the sum
system, is a unitary map for all
$t, a, b \in (0,\infty)$. (This would imply that the semigroup
$\{S_t\}$ can be extended as a semigroup of isometries on
$G_{(0,\infty)}$.) We denote by $A_{s,t}$ the map between $G_{(0,s)}
\oplus G_{(s, s+t)} \rightarrow G_{(0, s+t)}$ defined by $x \oplus y
\mapsto x +
y.$


\begin{Definition} We call a family
$\{x_t\}_{t \in (0,\infty)}$ such that
$x_t \in G_{(0,t)} , ~ \forall t \in (0,  \infty)$, as a real
additive unit for the sum system $(G_{(a,b)}, S_t)$, if

\medskip \noindent
(i) The map $ t \mapsto \langle x_t, x\rangle$ is a measurable map for
any $ x \in G_{(0,\infty)}$.

\medskip \noindent
(ii)
$x_s
+ S_s x_t = x_{s+t}, ~~ \forall s, t,
\in (0, \infty),$
(i. e.) $A_{s,t}(x_s \oplus S_s x_t)=x_{s+t}$.

Similarly we call a family
$\{y_t\}_{t \in (0,\infty)}$ such that
$y_t \in G_{(0,t)} , ~ \forall t \in (0,  \infty)$, as an imaginary
additive unit, for the sum system $(G_{(a,b)}, S_t)$,  if

\medskip \noindent
(i) The map $ t \mapsto \langle y_t, y\rangle$ is a measurable map for
any $ y \in G_{(0,\infty)}$.

\medskip \noindent
(ii) $\{y_t\}$ satisfies $(A_{s,t}^*)^{-1}(y_s \oplus S_s y_t) = 
y_{s+t},
~~
\forall s, t,
	\in (0, \infty).$
\end{Definition}

We denote by $R\A\U$ and $I\A\U$, the set of all real and imaginary
additive units respectively. For any given real(resptly. imaginary)
additive unit $\{x_t\}$ (resptly. $\{y_t\}$), we denote $x_{s,t} =
S_s(x_{t-s})
\in G_{(s, t)}$ (resptly. $y_{s,t} = S_s(y_{t-s})
\in G_{(s, t)}$).

We also define for an imaginary additive unit $\{y_t\}$, $$
 y^\prime_{s, s_1,s_2}= (A^*)^{-1}(0 \oplus y_{s_1, s_2} \oplus 0),
~\mbox{for any}~(s_1, s_2) \subset (0,s),$$ where $A:G_{(0,s_1)} \oplus
G_{(s_1,s_2)} \oplus G_{(s_2, s)} \mapsto G_{(0,s)}$, given by 
$x\oplus y \oplus z \mapsto
x + y+z$. It is easy to check that $y^\prime_{s, s_1,s_2} \in
G_{(0,s_1) \cup (s_2,
s)}^{\perp}.$ To simplify notation we denote 
$y^\prime _{1, s_1, s_2}$  by $y^\prime_{s_1, s_2}$, 
and
$y^\prime _{1,0, t}$ by
 just $y_t^\prime$. Finally note that $$
x_s+
x_{s,
s+t}= x_{s+t},~~
y_s^\prime +
y_{s,
s+t}^\prime = y_{s+t}^\prime.$$

\begin{Definition} A sum system $(G_{(a,b)},S_t)$ is called as a
divisible sum system if
the additive units exists and
generate the sum system,
(i. e.)
$$G_{(0,s)} =  \overline{span[x_{s_1,s_2}: s_1,s_2 \in (0,s), \{x_t\} \in
R\A\U]}$$ and $$G_{(0,s)} =  \overline{span[ y^\prime_{s, s_1,s_2}:
s_1,s_2 \in
(0,s), \{y_t\} \in I\A\U]}.$$
\end{Definition}

\bigskip

\begin{Proposition}\label{h}
\medskip \noindent (i) The collection of all real (and also imaginary)
aditive units forms a real vector space, with usual addition and scalar
multplication, $$\{x^1_t\} + \{x^2_t\} = \{x^1_t + x^2_t\};~~\lambda
\{x_t\}=\{\lambda x_t\}.$$

\medskip \noindent (ii) If $\{x_t\} \in R\A\U$ and $\{y_t\}
\in I\A\U$, then $$\langle x_t, y_t\rangle = \langle x_1, y_1\rangle t
~~\forall~ t \in (0, \infty).$$ In general for any two intervals $(s_1,
s_2),~(t_1, t_2) \subset (0, \infty)$, it is true that
\begin{eqnarray} \langle x_{s_1, s_2}, y^\prime_{t_1,
t_2}\rangle =
\langle x_1, y_1 \rangle~\ell
((s_1, s_2) \cap (t_1, t_2)),\end{eqnarray}
where $\ell$ is the Lebesgue
measure on
$\R$.

\medskip \noindent (ii) If a single real additive unit(and also an
imaginary additive unit) generates the sum system then the additive units
are determined uniquely up to a scalar.
\end{Proposition}

{\em Proof:} (i) Clear

\medskip \noindent (ii)Given any $\{x_t\} \in R\A\U$ and $\{y_t \} \in
I\A\U$, consider
the function $h_{x,y}(t)=\langle x_t,
y_t\rangle$. First we notice that $h_{x,y}$ is a real valued
measurable function. It may be proven as follows. We know that the map $t
\mapsto \langle x_t, x\rangle$ (also $t
\mapsto \langle y_t, x\rangle$) is measurable for any $ x \in
G_{(0,\infty)}$. Then $\|x\|= \sup_n \langle x_t, x_n\rangle$, for some
countable set $\{x_n\}$, due to the separability of the Hilbert space.
Hence
we
conclude that the function $t \mapsto
\|x_t\|$ is measurable. Similarly we conclude that the function $t \mapsto
\|y_t\|$ is also measurable. Now using the relation $$\langle x_t,
y_t\rangle = \frac{1}{4}(\|x_t + y_t\|^2 - \|x_t - y_t\|^2) ,$$ we can
conclude that the function $h_{x,y}(t)$ is measurable.

We also notice that $$h_{x,y}(s+t)=
\langle A_{s,t}(x_s \oplus
S_s(x_t)),
(A_{s,t}^*)^{-1}(y_s \oplus S_s y_t)\rangle = h_{x,y}(s) + h_{x,y}(t).$$
Therefore we conclude that
$h_{x,y}(t)=h_{x,y}(1) t$.
Now it is an easy verification to see that for
any two intervals $(s_1, s_2), (t_1, t_2) \subset (0, t)$
we have that
\begin{eqnarray}\label{xtyt} \langle x_{s_1, s_2}, y^\prime_{t_1,
t_2}\rangle =
h_{x,y}(1) ~\ell
((s_1, s_2) \cap (t_1, t_2)),\end{eqnarray} where $\ell$ is the Lebesgue
measure on
$\R$.

\medskip \noindent
(iii) Clear from (ii). \qed

\bigskip

\begin{Remark} (i) If the product system
is exponential, that
is the sum system $(G_{(a,b)}, S_t)$ is $(L^2((a,b), K), S_t)$, with the
standard shift $S_t$,
then $x_t = y_t =\xi 1_{(0,t)}$, for any $ \xi \in K$, exhausts all the
real and imaginary additive units, and the sum system is divisible.

\medskip\noindent (ii) The dimension of the vector space of additive
real (resptly. imaginary) units  may be defined as an index of the sum system, and it is clearly
an invariant for the sum system. In the case when the sum system gives
rise
to a type I product system it is a complete invariant. But in general it
is not, as all
examples of Tsirelson are of index $1$ and they are mutually
non-isomorphic.
\end {Remark}

We
prove in the next proposition that all
examples of
Tsirelson are divisible.

\begin{Proposition}\label{L2}
Let $(G_{(a,b)},S_t)$ be a sum system, and suppose that
$G_{(a,b)}$ is the completion of $L^2(a,b)$ with respect to some
inner product, such that $S_t$
the
canonical shift becomes an isometry. Then

\medskip
\noindent
(i) $x_t = 1_{0,t}$ is a real additive unit.

\medskip
\noindent
(ii) The non-zero imaginary additive unit exists (which is unique 
up to a scalar, if it
exists) if and only if the linear functional $f \mapsto
\int f dt$ is continuous on the dense subspace $L^2(a,b) \subset
G_{(a,b)}$.
\end{Proposition}

{\em Proof:}
As we have assumed the map $A_{s,t}:G_{(0,s)}~ \oplus ~
G_{(s,t)} \rightarrow
G_{(0, s+t)}$ to be $x \oplus y \mapsto  x +
y$, it is clear that $x_t = 1_{0,t}$ is a real additive unit and also it
generates the sum system.

To prove (ii), suppose a non-zero imaginary additive unit $\{y_t\}$
exists for the
sum system, then the relation \ref{xtyt}
(as $h(1) \neq 0$, and by choosing a real
multiple of
$y_t$ if needed) can be written, using (i), as $$\langle f, y_t\rangle =
\int_0^t f
dt,$$ for any simple function $f \in L^2(0,t) \subset G_{(0,t)}$.
Now it follows that
then the linear functional $f \mapsto
\int f$ is continuous on the dense subspace $L^2((a,b)) \subset
G_{(a,b)}$.
Suppose if we assume that the linear
functional $f \mapsto
\int f$ is continuous on the dense subspace $L^2((a,b)) \subset
G_{(a,b)}$, then we can choose $y_t \in G_{(0,t)}$ satisfying relation
\ref{xtyt} with $h(1)=1$. Now it is an easy verification to check that,
we have for $s,
t \in (0,\infty)$ and $s_1, s_2 \in
(0,s+t)$ that $$\langle x_{s_1, s_2},
(A^*)^{-1}(y_s \oplus S_s y_t)\rangle = \langle x_{s_1,s_2},
y_{s+t}\rangle.$$ As the set $\{x_{s,t}: s,t \in (0,s+t)\}$ is total in
$G_{(0,s+t)}$, we conclude that $\{y_t\}$ is an imaginary additive unit.
\qed

\bigskip

\begin{Corollary}
All examples of Tsirelson (Example \ref{Tsirelson}) are divisible.
\end{Corollary}

{\em Proof:}
To prove that sum systems of Example \ref{Tsirelson}) are
divisible,
we basically need to prove the existense of the the imaginary
additive unit,
(i.e.) it is enough to
prove that $f \mapsto \int f$ is continuous with respect to
the scalar product (\ref{sclrprdt}) for $f \in L^2(0,t)$.
That is we want a $ g \in L^2(0,t)$ such
that $g \star B = 1_{(0,t)}$, so
that $\int f g \star B=\int f$.
By taking
Fourier transform we basically
need a $\hat{g} \in
\ell^2(\Z)$, such that $\hat{g}\hat{B}= \hat{1}_{(0,t)}$, that is
we need to verify $\frac{e^{int}-1}{in\hat{B}} \in \ell^2(\Z)$.
But we have that $\hat{B}$ never vanishes and $\hat{B}(n) \sim
\frac{C}{\ln^{\alpha -1}|n|}$ for $n
\rightarrow \pm \infty$ (see \cite{Ts}, lemma 9.5, page 41), and the
series $\sum_{n \in \Z} \frac{\ln^{2\alpha -2}|n|}{n^2}$ is
convergent.
\qed

\bigskip

Now we prove that the product system arising from a divisible sum system
is always symmetric.

\begin{Proposition}
Suppose $(H_t, U_{s,t})$ be a product system
constructed out of a divisible sum system $(G_{(a,b)}, S_t)$, then $(H_t, U_{s,t})$ is
a symmetric product system.
\end{Proposition}

{\em Proof:} It is enough if we prove that the sum system is anti-isomorphic to itself.
Let $\{ \{x^i_t\}:i \in I\}$ be a spanning collection of real additive units for the sum
system.
Define $T_t: G_{(0,t)} \mapsto G_{(0,t)}$, by $T_t(x^i_{s_1, s_2})= x^i_{t-s_2, t-s_1},$ for
$(s_1, s_2) \subset (0,t), ~i \in I.$
Clearly $$\|T_t(x^i_{s_1, s_2})\|=\|x^i_{t-s_2, t-s_1}\| = \|S_{t-s_2}(x^i_{s_2
-s_1})\|=\|x^i_{s_1,
s_2}\|,$$ as we have assumed that the shift map to be isometric. So $T_t$ is
an isometry on a total set, and it is also bijective on this total set. Hence the map
$T_t$ extends
to a
unitary operator on $G_t$. It is easy to check that this map provides the required
anti-isomorphism. \qed

\bigskip

Next we prove a theorem which asserts only type $I$
and type
$III$ product systems can be constructed from a divisible sum
system.

\begin{Theorem}\label{typeI}
Let $(H_t, U_{s,t})$ be a product system
constructed out
of a divisible sum system $(G_{(a,b)}, S_t)$. If 
$(H_t, U_{s,t})$ has a unit then it is a  type $I$ product system.
\end{Theorem}

{\em Proof:} We assume that a unit $u(t) \in H_t =
\Gamma_s(\overline{G_{(0,t)}})$ exists for the product system, and prove
that
the product system is divisible.

Let $z_t \in \overline{G_{(0,t)}}$
be such that $z_t = c_1 x_t + ic_2  y_t$, where $\{x_t\} \in R\A\U,
\{y_t\} \in I\A\U$
and $c_1, c_2$ are real
scalars. Then clearly
it holds that
$S_{A_{s,t}}(z_s \oplus S_s z_t)=z_{s+t}$. So we
have $$U_{s,t}(W(z_s)\otimes W(z_t))U_{s,t}^*=\Gamma(A_{s,t})W(z_s
\oplus S_s z_t) \Gamma(A_{s,t})^*=W(z_{s+t}).$$ This basically
shows that
the family of unitaries $W(z_t) \in \B(H_t)$, is an automorphism
for the product system. As any automorphism of a product system preserves
units, we
conclude that the family of vectors $W(z_t) u_t \in H_t$ is also a
unit for the product
system $(H_t)$.

Fix a $t \in (0, \infty)$. The definition of
divisibility asserts that the set of all vectors of the form
$\sum_{j=1}^n c_j
x^j_{s_{j-1},s_j} + i c^\prime_j  {{y^\prime}^j_{t, s_{j-1},s_j}}$,
where $c_j,
c_j^\prime$
varying over real numbers, $s_0 =0 < s_1<s_2 \cdots< s_n=t$, and
$\{x^j_t\}$ and $\{y^j_t\}$ varying over all real and imaginary units
respectively, is dense
in
$\overline{G_{(0,t)}}$.

If we denote the unit $W(c x_t + ic^\prime  y_t)u_t$  by $\{v_{c, c^\prime}
(t)\}$, then
the image of $\otimes_{j=1}^n v_{c_j, c^\prime_j}(s_i-s_{i-1}),$ under the
canonical unitary of the product system is $W(\sum_{j=1}^n c_j
x_{s_{j-1},s_j} + i c^\prime_j  y^\prime_{t, s_{j-1},s_j})u_t$. So 
we
conclude that the units
generate the subspace $$\overline{span}[W(x)u_t: x \in \overline{G_{(0,t)}}].$$
But this subspace is whole of $\Gamma_s(\overline{G_{(0,t)}})$, as the
Weyl representation is
irreducible. Hence the product system is divisible, i.e. of type $I$.
\qed

\bigskip

For any
elementary set $E = \sqcup_{i=1}^n (s_i, s_{i+1}) \subset (0,1)$, we
define
$$x_E =\sum_{i=1}^n x_{s_i, s_{i+1}}, ~ y^\prime_E= \sum_{i=1}^n
y^\prime_{s_i, s_{i+1}} \in G_{(0,1)}~~~ \mbox{for}~\{x_t\}\in
R\A\U,~\{y_t\}
\in I\A\U.$$

The following theorem provides a necessary condition for the
product system arising from divisible sum system to be of type $I$. By
the next theorem the sufficient condition for the product system to
be of type $III$ is to violate this condition.

\begin{Theorem}\label{unitless} Let $(G_{(a,b)}, S_t)$ be a divisible sum
system, giving rise to a type $I$ product system. Then for any sequence
of elementary sets $E_n$
satisfying
$\liminf G_{E_n}= G_{(0,1)},~ \mbox{it also holds
that}~
\limsup G_{E_n^c}=\{0\}$
\end{Theorem}

{\em Proof:} Let $(H_t, U_{s,t})$ be the product system given by the sum
system $(G_{(a,b)}, S_t)$. As it is of type $I$ (see \cite{Arv}), it is
isomorphic to an exponential product system $(H_t^\prime,
U^\prime_{s,t})$, given by the sum system $(G^\prime_{(a,b)},
S^\prime_t)$,
where $G_{(a,b)} = L^2((a,b),K)$ for some separable Hilbert space $K$,
and $S^\prime_t$ is the canonical shift.
We denote by $(V_t)_{t \in (0, \infty)}$
a family of unitary maps implementing the isomorphism between the product
systems $H_t$ and $H_t^\prime$.

First let us note that the condition, $\liminf G^\prime_{E_n}=
G_{(0,1)}$ forcing $\limsup G^\prime_{E_n^c}=\{0\}$, is satisfied by the
sum
system $(G^\prime_{(a,b)},
S^\prime_t)$. This follows first by noticing that
$G^\prime_{E^c}={G^\prime_E}^{\perp}$ for
any
elementary set $E$, and then by
using lemma \ref{infsup}.

Now we claim that the set $\{y^\prime_B: B \in {\cal F}^e\}
\subset G_{(0,1)},$ is bounded for any imaginary additive unit $y$.
Suppose not, for each positive
integer $n$, choose an elementary set $B_n \subset (0,\frac{1}{n})$ such
that
$\|y^\prime_{B_n}\| > n$. If this is not possible for some $n$, that is
for each
elementary set $B \in (0,\frac{1}{n})$, $\|y^\prime_B\| \leq n$, then by
shifting
the
$y^\prime_B$s by the unitary operator $S_{\frac{k}{n}}, ~k=1, 2, \cdots
n$, and
by using the
triangle inequality we
may conclude that $\|y^\prime_E\| \leq n^2$, for any elementary set in $ E
\subset (0,1)$. But this means that the set $\{y^\prime_B: B \in {\cal
F}^e_{[0,1]}\}$ is
bounded. So we can indeed choose such $B_n \subset (0,\frac{1}{n})$ such
that $\|y^\prime_{B_n}\| > n$.

Now we know that $W(y_t)$ is an
automorphism for the product system $(H_t)$, hence $V_tW(y_t)V_t^*$ 
is an 
automorphism of the product system $(H_t^\prime)$. By the result in
section
8 of \cite{Arv} we can conclude that $$V_tW(y_t)V_t^* = e^{i\lambda
t}W(\xi 1_{(0,t)} )Exp(U^t), ~\mbox{where}~ \lambda \in \R, ~\xi \in 
K, ~U \in
\U(K),$$ and $U^t( \eta 1_A)= (U\eta)1_A$ for any $\eta \in K$ and $A
\subset
(0,t)$.  It is easy to verify that $$V_1 W(y^\prime_{B_n}) V_1^*=
e^{i\ell(E_n)}W(\xi
1_{B_n}) Exp(U^{E_n}),$$ where $U^{E_n}$ is the unitary operator 
defined
by $U^{E_n} (\eta 1_{A})= U\eta 1_A$ if $ A \subset E_n$, and 
$U^{E_n}  
(\eta 1_{A^\prime})= \eta 1_A^\prime$ if $ A^\prime \subset E_n^c$, 
$\eta
\in K$. Clearly the above
sequence converges strongly to the identity operator.
But the
sequence $W(y^\prime_{B_n})$ can not be a strongly convergent 
sequence,  
for the following reason. Suppose $W(y^\prime_{B_n})$ converges 
strongly,
then by applying on the vacuum vector, we
first conclude that the sequence 
$e^{-\frac{\|y^\prime_{B_n}\|^2}{2}}e(y^\prime_{B_n})$
converges. But the projection of this sequence on to the k-th
particle subspace
converges to $0$,
for each $k \in \N$, as
$e^{-\frac{\|y^\prime_{B_n}\|^2}{2}} 
(\frac{\|y^\prime_{B_n}\|^k}{\sqrt{n!}})$
converges to $0$, for each $k \in \N$. Hence the sequence
$e^{-\frac{\|y^\prime_{B_n}\|^2}{2}}e(y^\prime_{B_n})$ should 
converge to $0$, 
but
this is
not possible as 
$e^{-\frac{\|y^\prime_{B_n}\|^2}{2}}\|e(y^\prime_{B_n})\| = 1$ 
for
all $n$. Hence we have proved our claim.

Now suppose $E_n$ be any sequence of elementary sets satisfying 
$\liminf
G_{E_n}= G_{(0,1)}$, we claim that $\ell(E_n^c)$ converges to
$0$,
as $n$ tends to $\infty$. Suppose $\ell(E_n^c)$ does not converge to 
$0$,
let $y^\prime_{n}=y^\prime_{E_{n}^c} \in G_{E_{n}}^{\perp}.$ Then 
$\|y_n^\prime\|$ is a bounded sequence which
does not converge weakly
to $0$, as
$$\langle y_n^\prime, x_1\rangle= \ell(E_{n}^c).$$
So we conclude that $\limsup G_{E_n}^\perp$ is not equal to
$\{0\}$. But, by lemma \ref{infsup}, this contradicts our assumption 
that
$\liminf G_{E_n}= G_{(0,1)}$. Hence we have proved our claim that
$\ell(E_n^c) \rightarrow 0$.

As any $f \in L^2((0,1), K)$ is the limit of $1_{E_n} f$, we have
that $$\liminf L^2(E_n, K)=
L^2((0,1), K).$$ But this also implies that $\limsup L^2(E_n^c, K)
= \{0\}.$ Now the theorem \ref{invariant} implies that $\limsup
G_{E_n^c} = \{0\}$ and the proof of the proposition is over.
\qed

\bigskip

\begin{Remark}\label{Tsirelunitless}
Tsirelson in his examples produces a sequence of elementary sets
$\{E_n\}$ such that
$\liminf G_{E_n}= G_{(0,1)},$ but the condition that
$\limsup G_{E_n^c}=\{0\}$ is violated. This once again proves that the examples of
Tsirelson are of type $III$.
\end{Remark}

\noindent{\bf Acknowledgements:} The first author is supported by
the Swarnajayanthi Fellowship from the Department of Science and
Technology (India). The second author is currently a JSPS 
postdoctoral fellow and  earlier, when this work got initiated, was
 supported by a postdoctoral
fellowship from IFCPAR grant No. IFC/2301-1/2001/1253 
at the Indian
Statistical Institute.

\end{document}